\documentclass[11pt,a4paper]{article}

\usepackage[english]{babel}
\usepackage[utf8]{inputenc}
\usepackage{amsmath}
\usepackage{amsfonts}
\usepackage{graphicx}
\usepackage{epstopdf}
\usepackage{enumitem}
\usepackage[font=footnotesize, labelfont=bf]{caption}
\usepackage[figbotcap, loose]{subfigure}

\usepackage{siunitx}
\usepackage{tabularx}
\usepackage{booktabs}
\usepackage{multirow}
\usepackage{hyperref}

\usepackage{color}
\usepackage{soul}

\usepackage[square,numbers,sort&compress]{natbib}

\def\pier #1{{\color{black}#1}}
\def\pier #1{#1}
\def\pg #1{{\color{black}#1}}
\def\pg #1{#1}
\def\pierbis #1{{\color{black}#1}}

\def\vec#1{\boldsymbol{#1}}
\usepackage{placeins}

\usepackage{chngcntr}
\counterwithin*{equation}{section}

\begin{document}

\begin{center}
\textbf{\Large Optimal control of cytotoxic and antiangiogenic therapies on prostate cancer growth}

\vspace{5mm}

{Pierluigi Colli}

{\footnotesize
{Dipartimento di Matematica, Universit\`a degli Studi di Pavia and IMATI-C.N.R., \\
Via Ferrata~5, 27100 Pavia, Italy\\
pierluigi.colli@unipv.it
}}

\vspace{2mm}

{Hector Gomez}

{\footnotesize
{School of Mechanical Engineering, Purdue University,\\
516 Northwestern Avenue, West Lafayette, IN 47907, USA\\
and\\
Weldon School of Biomedical Engineering, Purdue University, \\
206 S. Martin Jischke Drive, West Lafayette, IN 47907, USA\\
and\\
Purdue Center for Cancer Research, Purdue University, \\
201 S. University Street, West Lafayette, IN 47907, USA\\
hectorgomez@purdue.edu}}

\vspace{2mm}

{Guillermo Lorenzo}

{\footnotesize
{
Oden Institute for Computational Engineering and Sciences,\\
The University of Texas at Austin,\\
201 E. 24th Street, Austin, TX 78712-1229, USA\\
and\\
Dipartimento di Ingegneria Civile e Architettura, \\
Universit\`{a} degli Studi di Pavia and IMATI-C.N.R.,\\
Via Ferrata 3, 27100 Pavia, Italy \\
guillermo.lorenzo@unipv.it}}

\vspace{2mm}

{Gabriela Marinoschi}

{\footnotesize
{``Gheorghe Mihoc-Caius Iacob'' Institute of Mathematical Statistics \\
and Applied Mathematics of the Romanian Academy,\\
Calea 13 Septembrie 13, 050711 Bucharest, Romania\\
gabriela.marinoschi@acad.ro}}

\vspace{2mm}

{Alessandro Reali}

{\footnotesize
{Dipartimento di Ingegneria Civile e Architettura, \\
Universit\`{a} degli Studi di Pavia and IMATI-C.N.R.,\\
Via Ferrata 3, 27100 Pavia, Italy \\
alessandro.reali@unipv.it}}

\vspace{2mm}

{Elisabetta Rocca}

{\footnotesize
{Dipartimento di Matematica, Universit\`a degli Studi di Pavia and IMATI-C.N.R.,\\
Via Ferrata~5, 27100 Pavia, Italy\\
elisabetta.rocca@unipv.it
}}

\vspace{5mm}

\end{center}

\begin{center}
{\textbf{Abstract}}
\end{center}
Prostate cancer can be lethal in advanced stages, for which chemotherapy may become the only viable therapeutic option. While there is no clear clinical management strategy fitting all patients, cytotoxic chemotherapy with docetaxel is currently regarded as the gold standard. However, tumors may regain activity after treatment conclusion and become resistant to docetaxel. This situation calls for new delivery strategies and drug compounds enabling an improved therapeutic outcome. Combination of docetaxel with antiangiogenic therapy has been considered a promising strategy. Bevacizumab is the most common antiangiogenic drug, but clinical studies have not revealed a clear benefit from its combination with docetaxel. Here, we capitalize on our prior work on mathematical modeling of prostate cancer growth subjected to combined cytotoxic and antiangiogenic therapies, and propose an optimal control framework to robustly compute the drug-independent cytotoxic and antiangiogenic effects enabling an optimal therapeutic control of tumor dynamics. We describe the formulation of the optimal control problem, for which we prove the existence of at least a solution and determine the necessary first order optimality conditions. We then present numerical algorithms based on isogeometric analysis to run a preliminary simulation study over a single cycle of combined therapy. Our results suggest that only cytotoxic chemotherapy is required to optimize therapeutic performance and we show that our framework can produce superior solutions to combined therapy with docetaxel and bevacizumab. We also illustrate how the optimal drug-na\"{i}ve cytotoxic effects computed in these simulations may be successfully leveraged to guide drug production and delivery strategies by running a nonlinear least-square fit of protocols involving docetaxel and a new design drug. In the future, we believe that our optimal control framework may contribute to accelerate experimental research on chemotherapeutic drugs for advanced PCa and ultimately provide a means to design and monitor its optimal delivery to patients.

\vspace{2mm}
{\small

\vspace{2mm}
}


\section{Introduction}

Prostate cancer (PCa) is the second most common cancer among men worldwide \cite{ferlay2015cancer}. Current estimates indicate that the annual number of new cases globally exceeds 1,000,000. PCa is also responsible for over 300,000 worldwide deaths yearly and a major economic burden for our health care systems. 

In most cases, PCa is an adenocarcinoma  --- a form of cancer that originates and develops in epithelial tissues. PCa is easier to treat in its early stages, especially before it spreads out of the prostate \cite{walsh2}. Although the disease rarely produces any symptoms before the tumor has advanced significantly, current screening protocols have helped identify PCa patients early \cite{Mottet2018,walsh2}. The most common screening method for PCa is the Prostate Specific Antigen (PSA) test --- a blood test that measures the serum level of gamma-seminoprotein which is a biomarker of the prostate activity. PSA-based screening is controversial in the urology community, primarily because it is prostate-specific, but not cancer-specific \cite{walsh2}. Despite the success of screening protocols, a significant number of patients eventually develops advanced PCa \cite{Buzzoni2015,walsh2}. These include men who are diagnosed late, and patients who initially receive radical treatment (e.g., surgery or radiation) or follow a conservative management protocol (e.g., active surveillance), but experience cancer progression.

The most common treatments for advanced PCa are hormonal therapy and chemotherapy, but the latter may become the only remaining alternative after PCa becomes insensitive to hormonal strategies \cite{Mottet2018,walsh2}.
Here, we will focus on chemotherapy, which normally relies on cytotoxic drugs.
These compounds aim at interrupting cell division and promoting cell death, hence contributing to reduce tumor progression.
Cytotoxic chemotherapy for advanced PCa usually relies on the periodic delivery of docetaxel every three weeks \cite{Mottet2018,walsh2,Mottet2018,Seruga2011,Eisenberger2012}.
However, treated tumors may resume growth months after the conclusion of therapy and  develop resistance to docetaxel, hence requiring the use of alternative cytotoxic drugs for a second round of treatment.
Angiogenesis has also been regarded as a promising therapeutic target for advanced PCa because it is known to play a pivotal role in the progression of advanced PCa and higher microvascular density correlates with poorer prognosis \cite{weidner1993tumor,Mucci2009,mehta2001independent}. 
Antiangiogenic drugs block the formation of new blood vessels cutting off the supply of critical nutrients to the tumor.
Bevacizumab is an extensively studied antiangiogenic drug that has also being investigated for the treatment of advanced PCa either as monotherapy or in combination with cytotoxic chemotherapy \cite{Antonarakis2012,Seruga2011,Small2012,mukherji2013angiogenesis}.
A clinical study showed that castration-resistant PCa patients did not benefit from treatment with bevacizumab alone \cite{reese2001phase}, but another study showed that 75\% of the patients treated with bevacizumab experienced a 50\% PSA decline when taking bevacizumab in combination with cytotoxic chemotherapy \cite{picus2011phase}.
However, multiple experimental and clinical studies have shown contradictory evidence about the effectiveness of antiangiogenic therapy for advanced PCa either as monotherapy or combined with a cytotoxic agent \cite{Antonarakis2012,Seruga2011,Small2012,mukherji2013angiogenesis,Kelly2012,Tannock2013,Petrylak2015,Smith2016,Sternberg2016}. 

Thus, developing optimal dosing strategies for chemotherapeutic drugs is an open problem of major medical interest, but it is particularly difficult when two or more drugs are combined because their joint effect may produce unexpected outcomes and toxicities. 
Recently, the use of mathematical models to describe the growth and treatment of cancer has been increasing our understanding of these pathologies and providing a means to monitor and forecast tumor growth \emph{via} patient-specific computer simulations, which can assist physicians in clinical decision-making \cite{Anderson2008,Yankeelov2013,Corwin2013,Byrne2006,BLM,Oden2010,Rocha2018,Yin2019}.
Several studies have focused on the effects of chemotherapeutic strategies for various types of tumors \cite{Bogdanska2017,Jarrett2018,Gallaher2018,Hinow2009,Kohandel2007,Powathil2007,IrurzunArana2020,Yin2019}.
In this context, optimal control theory provides a robust framework to obtain optimal drug protocols according to a given mathematical model of tumor growth and treatment \cite{Jarrett2020,Benzekry2013,CRW,GLR,CGRS1,Colli2019,Ledzewicz2007,Fister2003,Leszczynski2019,Wang2016}.
We argue that these quantitative approaches can play a significant role in the design of optimal drug regimens, providing the best dosing, timing, and even drug pharmacodynamics to treat each patient's tumor \cite{Jarrett2020,IrurzunArana2020,Luepfert2005,Yin2019}.
Hence, optimal drug protocols can help us reduce the total amount of drug needed to achieve a target tumor size reduction, which is very important considering the harsh side effects of cytotoxic and antiangiogenic drugs \cite{Mottet2018,Seruga2011,Eisenberger2012,Antonarakis2012,Seruga2011,Small2012,mukherji2013angiogenesis}. Because antiangiogenic and cytotoxic drugs attack the tumor using different mechanisms, mathematical modeling can also help us determine the most effective drug dose combination.
Additionally, optimal pharmacodynamics derived from optimal control problems could ultimately be exploited to guide the design of new drug compounds \cite{Luepfert2005,Iyengar2012,Hussain2019,Shi2017}.

Here, we propose an optimal control framework that relies on a recently proposed model of PCa growth accounting for the effect of cytotoxic and antiangiogenic drugs \cite{CGLMRR-1}.  Our approach aims at finding the cytotoxic and antiangiogenic drug effects minimizing an array of common quantities of interest in clinical and experimental settings, such as tumor volume and serum PSA.
We demonstrate the existence of unique solutions and derive the optimality conditions for our optimal control problem. Then, we propose a set of numerical algorithms to solve this optimal control problem and we use them to run a simulation study. Our computational methods rely on isogeometric analysis (IGA), a recent generalization of the finite element method with superior approximation properties granted by the use of splines as basis functions \cite{hughes2005isogeometric}.
Our simulation study results in optimal therapeutic solutions that only feature cytotoxic effects, and hence suggests that antiangiogenic therapies may not be optimal for treating PCa.
Because docetaxel and bevacizumab are reference drugs in therapeutic investigations of advanced PCa, we used their properties to build the initial guesses and admissible space of controls.
However, our optimal control problem is drug-independent, i.e., it provides optimal cytotoxic and antiangiogenic drug effects without accounting for any specific drug delivery plan. Thus, in this work, we also capitalize on this advantageous property and propose a two-stage methodology to derive optimal drug compounds and delivery plans matching the optimal drug effects targeting a specific tumor: first, we calculate the optimal drug-na\"{i}ve optimal effects \textit{via} computer simulation of our optimal control problem, and, second, we run nonlinear least-square fits of an array of elective drug delivery plans and compounds. Our results show that our approach can render interesting optimal drug protocols matching the optimal drug effects to treat the target tumor.

The paper is organized as follows. Section~\ref{mathmodel} briefly describes our mathematical model for PCa growth featuring cytotoxic and antiangiogenic drug effects, as well as the formulation of the optimal control problem. Sections~\ref{exist} and \ref{optcon} present the existence results and the optimality conditions for the optimal control problem, respectively. Section~\ref{numerics} introduces the numerical algorithms that we propose in order to solve it, and we use these computational methods in a simulation study presented in Section~\ref{simstudy}. Then, we use the resulting optimal drug effects to illustrate how they can guide drug design and the conception of new drug delivery strategies in Section~\ref{drugprotocols}. Finally, in Section~\ref{disc} we discuss the work presented herein, draw conclusions, and address future avenues of research.


\section{Mathematical model}\label{mathmodel}

\subsection{Prostate cancer model with cytotoxic and antiangiogenic therapies}\label{pcamodel}

We leverage a recently proposed mathematical model of PCa growth \cite{CGLMRR-1} that includes the effect of cytotoxic and antiangiogenic drugs. Reference \cite{CGLMRR-1} includes a detailed description of the biological mechanisms considered in the model, and a biological interpretation of its parameters. Let $\Omega\subset\mathbb{R}^d$ be an open bounded set and $d$ the number of spatial dimensions. Let $\phi_0$, $\sigma_0$, and $p_0$ be sufficiently smooth given functions defined in $\Omega$. The set $\Omega$ has a sufficiently smooth boundary $\partial \Omega$ and $T$ is a finite time. By defining the space/time sets $Q_{T}:=(0,T)\times \Omega$ and $\Sigma _{T}:=(0,T)\times \partial \Omega$, the model can be formulated as 
\begin{alignat}{2}
&\phi _{t}-\lambda \Delta \phi +F^{\prime }(\phi )-m(\sigma )h^{\prime }(\phi)=-Uh^{\prime }(\phi ) &\quad\mbox{in }&Q_{T} ,  \label{1} \\
&\sigma _{t}-\eta \Delta \sigma +\gamma _{h}\sigma +\gamma _{ch}\sigma \phi = %
\pg{S_h (1-\phi) + (S_c - S)\phi} &\quad\mbox{in }&Q_{T},  \label{2} \\
&p_{t}-D\Delta p+\gamma _{p}p=\alpha _{h}+\alpha _{ch}\phi &\quad\mbox{in }&Q_{T},
\label{3} \\
&\phi =0,\mbox{ }\frac{\partial \sigma }{\partial \nu }=\frac{\partial p}{\partial\nu }=0 & \quad\mbox{in }&\Sigma _{T},  \label{4}\\
&\phi(0,x) =\phi_{0}(x),\mbox{ }\sigma(0,x)=\sigma_{0}(x),\mbox{ }p(0,x)=p_{0}(x) &\quad \mbox{in }&\Omega ,  \label{5}
\end{alignat}
where $\phi$ is a phase-field that identifies the spatial location of the tumor; $\sigma$ is the concentration of a vital nutrient; and $p$ is the tissue PSA concentration. The serum PSA $P_s$ commonly used in clinical practice can be obtained as $P_s=\int_\Omega p\, dx$ \cite{CGLMRR-1,lorenzo2016tissue}. The differential operators in the equations are defined as follows: the subscript $t$ indicates partial differentiation with respect to time; $\Delta$ denotes the Laplace operator with respect to the space variables; $F^\prime$ denotes the derivative of $F$ and $\frac{\partial }{\partial \nu }$ represents the outward normal derivative to $\partial \Omega $. In Eq.~\eqref{1}, $F(\phi )=M\phi ^{2}(1-\phi )^{2}$ and $h(\phi)=M\phi ^{2}(3-2\phi)$, where $M$ is a positive constant. We will also use the notation
\begin{equation}
\gamma _{ch}:=\gamma _{c}-\gamma _{h},\mbox{ }\ \alpha _{ch}:=\alpha_{c}-\alpha _{h}, \ \ \pg{S_{ch}:= S_c - S_h.}
\end{equation}
The model parameters $\lambda ,$ $\eta ,$ $D,$ $\gamma_p$, $\gamma _{c},$ $\gamma _{h},$ $\alpha _{c},$ $\alpha _{h}$, $S_c$ and $S_h$ are positive constants with biological significance; see Ref.~\cite{CGLMRR-1}. For the solution of the initial boundary value problem \eqref{1}--\eqref{5}, the functions $U(t,x)$ and $S(t,x)$ are considered to be given. 
$U(t,x)$ models the inhibiting effect of a cytotoxic drug on tumor dynamics, while $S(t,x)$ represents the reduction of nutrient supply produced by an antiangiogenic drug. The function $m(\sigma)$ is defined as
\begin{equation}\label{msigma}
m(\sigma )=m_{ref}\left( \frac{\rho +A}{2}+\frac{\rho -A}{\pi }\arctan\left( \frac{\sigma -\sigma _{l}}{\sigma _{r}}\right) \right)
\end{equation}
where $m_{ref}$ is a positive constant, while $\rho$ and $A$ are constants that determine the rate of cell proliferation and apoptosis. We define $\rho=K_\rho/\overline{K}_\rho$ and $A=-K_A/\overline{K}_A$. Here, $K_\rho$ and $K_A$ are, respectively, the proliferation and apoptosis rates of tumor cells while $\overline{K}_\rho$ and $\overline{K}_A$ are scaling positive constants. The positive constants $\sigma_l$ and $\sigma_r$ are, respectively, a reference and a threshold value for the nutrient concentration; see Ref.~\cite{CGLMRR-1}.

\subsection{Optimal control problem}\label{optcontrol}
Our optimal control problem consists in finding the functions $U(t,x)$ and $S(t,x)$ that provide the optimal cytotoxic and antiangiogenic effects to treat a certain PCa case. Therefore, our optimal control problem can be stated as
\begin{equation*}
\mbox{Minimize }\{J(U,S);\mbox{ }(U,S)\in \mathcal{U}_{ad}\},\tag*{($P$)}
\end{equation*}%
\noindent subject to Eq.~\eqref{1}--\eqref{5}, where
\begin{alignat}{2}
J(U,S) &=\frac{k_{1}}{2}\int_{Q_{T}}(\phi (t,x)-\phi _{Q})^{2}dxdt+\frac{k_{2}}{2}\int_{\Omega }(\phi (T,x)-\phi _{\Omega })^{2}dx  +k_{3}\int_{\Omega }\phi (T,x)dx \nonumber \\
&+\frac{k_{4}}{2}\int_{0}^{T}\left[ \left(\int_{\Omega }p(t,x)dx-p_{\Omega }(t)\right) ^{+}\right] ^{2}dt+k_{5}\int_{\Omega }p(T,x)dx  \notag \\
&+\frac{k_{6}}{2}\int_{Q_{T}}U^{2}(t,x)dxdt+\frac{k_{7}}{2}\int_{Q_{T}}S^{2}(t,x)dxdt  \label{J}
\end{alignat}%
\noindent and 
\begin{equation}
\mathcal{U}_{ad}=\{(U,S)\in L^{\infty }(Q_{T})\times L^{\infty }(Q_{T});\mbox{ }0\leq U\leq U_{\max },\mbox{ }0\leq S\leq S_{\max },\mbox{ a.e. in }%
Q_{T}\}.  \label{Uad}
\end{equation}

The quantities $U_{\max }$ and $S_{\max }$ are two threshold positive values such that $S_{\max} \leq S_c$. In Eq.~\eqref{J}, $r^{+}$ represents the positive part of $r$; the functions $\phi _{Q}, $ $\phi _{\Omega }$ are prescribed targets for the tumor phase field in $Q_{T}$ and in $\Omega$ at the final time, respectively, while $ p_{\Omega }$ represents an upper target function for the spatial mean value of the tissue PSA (i.e., an upper target for the serum PSA). We assume that 
\begin{equation}
\phi _{Q}\in L^{2}(Q_{T}),\mbox{ }\phi _{\Omega }\in L^{2}(\Omega ),\mbox{ }%
p_{\Omega }\in L^{2}(0,T).  \label{5-0}
\end{equation}

Notice that we have not imposed any specific drug delivery plan or pharmacodynamic properties for $U(t,x)$ and $S(t,x)$ in our mathematical model and the optimal control problem herein, hence making our approach drug-independent and solely relying on the type of therapy itself. In this context, $\mathcal{U}_{ad}$ aims at defining feasible limits for the therapeutic effects that can be realistically achieved with common drugs.
The optimal choice of delivery dates, doses and even pharmacodynamic properties can be made at a second stage, aiming at reproducing the optimal drug-na\"{i}ve effects obtained from the optimal control problem described in this section (see Section~\ref{drugprotocols}).

The coefficients $k_{i},$ $i=1,...,7,$ in Eq.~\eqref{J} are nonnegative constants, such that there exists at least one $i$ with $k_{i}>0$. In the numerical simulations presented in Section~\ref{simstudy}, we will focus on particular examples of the functional $J$, where only a few of the $k_i$'s are simultaneously nonzero. For example, the first three terms in $J$, weighted by $k_1$, $k_2$, and $k_3$, respectively, provide a slightly different control of the tumor extent. In our numerical simulations, not all three constants $k_1$, $k_2$, and $k_3$ will be nonzero simultaneously. The same argument can be applied to the constants $k_4$ and $k_5$, which provide control of the PSA. The terms weighted by $k_6$ and $k_7$ control the relative participation of cytotoxic and antiangiogenic effects, which may ultimately control the total amount of drugs delivered (see Ref.~\cite{CGLMRR-1} and Sections~\ref{simstudy}--\ref{drugprotocols}). 
The control of the total amount of drug may be especially important because some therapeutic compounds have very strong side effects. However, specific dose-dependent toxicities can be accounted for during the design of optimal drug protocols aiming at reproducing the optimal drug-na\"{i}ve effects obtained from our optimal control problem (see Sections~\ref{simstudy}, \ref{drugprotocols}, and \ref{disc}).

\section{Existence results}\label{exist}

\setcounter{equation}{0}

The functional framework involves the space $H=L^{2}(\Omega ),$ which is
identified with its dual space $H^{\prime }\equiv H,$ and the Sobolev spaces 
\pg{$H^{1}(\Omega )$ and $H_{0}^{1}(\Omega )$,
the last one containing the elements of  $H^{1}(\Omega )$ with null trace on the boundary $\partial \Omega$. We set}
\begin{equation*}
V_{0}=H_{0}^{1}(\Omega ),\mbox{ }V_{0}^{\prime }=(H_{0}^{1}(\Omega
))^{\prime }=:H^{-1}(\Omega ),\mbox{ \ }V=H^{1}(\Omega ),\mbox{ }V^{\prime
}=(H^{1}(\Omega ))^{\prime },\mbox{ }
\end{equation*}
\pg{as well as} 
\begin{equation*}
W_{0}=H^{2}(\Omega )\cap H_{0}^{1}(\Omega ),\mbox{ \ }W=\left\{ y\in
H^{2}(\Omega );\mbox{ }\frac{\partial y}{\partial \nu }=0\mbox{ on }\partial
\Omega \right\} ,
\end{equation*}%
with the dense and compact injections $W_{0}\subset V_{0}\subset H\subset
V_{0}^{\prime }$ and $W\subset V\subset H\subset V^{\prime }.$ 
\pg{For the norms in these spaces we will use the notation $\| \, \cdot \, \|_B$, where 
$B$ is the space we are considering.}

\pg{For $r\in \lbrack 1,\infty ]$ and $z\in L^{r}(\Omega )$ or $z\in
L^{r}(Q_{T}) $ we simply denote the norms of $z$ by~$\left\Vert z\right\Vert_{r}.$}

\pg{We will also 
make use of spaces of functions that depend on time with values in a Banach space $B$. Namely, for 
 $r\in \lbrack 1,\infty ]$ 
we consider the space $L^r (0,T;B)$ of measurable functions $t\mapsto z(t)$ such that 
$t \mapsto \| z(t)\|_B^r $ is integrable on $(0,T)$ (or essentially bounded  if $r=\infty$) 
and the space $C([0,T]; B) $ of continuous functions from $[0,T]$ to $B$. Perhaps, it is important to 
note  that $L^2 (0,T;H)$ is a space completely isomorphic to $L^2 (Q_T)$. Moreover, for 
 $r\in \lbrack 1,\infty ]$,  $W^{1,r} (0,T;B)$
will denote the space of functions $t\mapsto z(t)$ such that both $z$ and its (weak) derivative $z_t$
belong to $L^r (0,T;B)$. We point out that  $W^{1,r} (0,T;B) \subset  C([0,T]; B) $
for all  $r\in \lbrack 1,\infty ]$.} \medskip

\textbf{Hypotheses.} According to the explanation 
\pg{already given, we
emphasize} that%
\begin{equation}
F\in C^{\infty }(\mathbb{R}),\mbox{ }h\in C^{\infty }(\mathbb{R}),\mbox{ }%
m\in C^{\infty }(\mathbb{R}),  \label{6}
\end{equation}
\begin{equation}
m,\mbox{ }m^{\prime }\mbox{ are Lipschitz continuous on }\mathbb{R},
\label{7}
\end{equation}
\begin{equation}
m,\mbox{ }m^{\prime },\mbox{ }m^{\prime \prime }\in L^{\infty }(\mathbb{R}).
\label{8}
\end{equation}

Let us introduce the spaces%
\begin{eqnarray}
X_{0} &=&W^{1,2}(0,T;V_{0}^{\prime })\cap C([0,T];H)\cap L^{2}(0,T;V_{0}),
\notag \\
X &=&W^{1,2}(0,T;V^{\prime })\cap C([0,T];H)\cap L^{2}(0,T;V),  \label{9}
\end{eqnarray}%
\begin{eqnarray}
\mathcal{X}_{0} &=&W^{1,2}(0,T;H)\cap C([0,T];V_{0})\cap L^{2}(0,T;W_{0}),
\notag \\
\mathcal{X} &=&W^{1,2}(0,T;H)\cap C([0,T];V)\cap L^{2}(0,T;W).  \label{10}
\end{eqnarray}

\medskip

\noindent \textbf{Definition 3.1.} A solution to system (\ref{1})-(\ref{5})
is a triplet $(\phi ,\sigma ,p),$ with $\phi \in X_{0}\cap L^{\infty
}(Q_{T}),$ $\sigma \in X,$ $p\in X,$ which satisfies 
\begin{eqnarray}
&&\int_{0}^{T}\left\langle \phi _{t}(t),\psi (t)\right\rangle
_{V_{0}^{\prime },V_{0}}dt+\int_{Q_{T}}\left\{ \lambda \nabla \phi \cdot
\nabla \psi +(F^{\prime }(\phi )-m(\sigma )h^{\prime }(\phi ))\psi \right\}
dxdt\;\mbox{ \ \ \ }  \notag \\
&=&-\int_{Q_{T}}Uh^{\prime }(\phi )\psi dxdt\mbox{, for all }\psi \in
(0,T;V_{0}),  \label{e9-1}
\end{eqnarray}%
\begin{eqnarray}
&&\int_{0}^{T}\left\langle \sigma _{t}(t),\psi (t)\right\rangle _{V^{\prime
},V}dt+\int_{Q_{T}}(\eta \nabla \sigma \cdot \nabla \psi +(\gamma _{h}\sigma
\psi +\gamma _{ch}\sigma \phi \psi )dxdt  \notag \\
&=&\int_{Q_{T}}\pg{(S_h + (S_{ch} -S) \phi )}\psi dxdt,\mbox{ for all }\psi
\in L^{2}(0,T;V),  \label{e9-2}
\end{eqnarray}%
\begin{eqnarray}
&&\int_{0}^{T}\left\langle p_{t}(t),\psi (t)\right\rangle _{V^{\prime
},V}dt+\int_{Q_{T}}(D\nabla \pierbis{p} \cdot \nabla \psi +\gamma _{p}p\psi )dxdt
\notag \\
&=&\int_{Q_{T}}(\alpha _{h}+\alpha _{ch}\phi )\psi dxdt,\mbox{ for all }\psi
\in L^{2}(0,T;V),  \label{e10}
\end{eqnarray}%
and 
\begin{equation}
(\varphi ,\sigma ,p)(0)=(\varphi _{0},\sigma _{0},p_{0})\quad%
\hbox{a.e. in
$\Omega$}.  \label{e10-0}
\end{equation}

\medskip

\noindent \textbf{Theorem 3.2. }\textit{Let }$(U,S)\in \mathcal{U}_{ad},$%
\textit{\ } \textit{\ }%
\begin{equation}
(\phi _{0},\sigma _{0},p_{0})\in H\times H\times H,  \label{e12}
\end{equation}%
\begin{equation}
0\leq \phi _{0}(x)\leq 1\mbox{ \textit{a.e.} }x\in \Omega .  \label{e12-0}
\end{equation}%
\textit{Then, system} (\ref{1})-(\ref{5}) \textit{has a unique solution }$%
(\phi ,\sigma ,p),$ \textit{with} 
\begin{equation}
0\leq \phi (t,x)\leq 1\mbox{ \textit{a.e.} }(t,x)\in Q_{T},  \label{e12-1}
\end{equation}%
\textit{which satisfies the estimate}%
\begin{eqnarray}
&&\left\Vert \phi \right\Vert^2_{X_{0}}+\left\Vert \sigma \right\Vert^2_{X}+\left\Vert p\right\Vert^2_{X}  \notag \\
&\leq &C(\left\Vert \phi _{0}\right\Vert _{H}^{2}+\left\Vert \sigma
_{0}\right\Vert _{H}^{2}+\left\Vert p_{0}\right\Vert _{H}^{2}+\left\Vert
U\right\Vert _{L^{2}(0,T;H)}^{2}+\left\Vert S\right\Vert
_{L^{2}(0,T;H)}^{2}+1).  \label{e12-00}
\end{eqnarray}%
\textit{If \ } 
\begin{equation}
(\sigma _{0},p_{0})\in L^{\infty }(\Omega )\times L^{\infty }(\Omega ),%
\mbox{
}\sigma _{0}(x)\geq 0,\mbox{ }p_{0}(x)\geq 0\mbox{ \textit{a.e.} }x\in
\Omega ,  \label{e13-0}
\end{equation}%
\textit{then we have }$(\sigma ,p)\in L^{\infty }(Q_{T})\times L^{\infty
}(Q_{T}),$ 
\begin{equation}
\sigma (t,x)\geq 0,\mbox{ }p(t,x)\geq 0\mbox{ \textit{a.e.} }(t,x)\in Q_{T},
\label{e13-2}
\end{equation}%
\textit{and }$(\sigma ,p)$\textit{\ fulfills the estimate} 
\begin{equation}
\left\Vert \sigma \right\Vert _{\infty }\leq C(\left\Vert \sigma
_{0}\right\Vert _{\infty }+\pg{1}),\mbox{ \ }\left\Vert p\right\Vert
_{\infty }\leq C(\left\Vert p_{0}\right\Vert _{\infty }+1).  \label{e13-3}
\end{equation}%
\textit{Moreover, the solution is continuous with respect to the data, that
is, for two solutions }$(\phi _{i},\sigma _{i},p_{i})$\textit{\
corresponding to }$(\phi _{0}^{i},\sigma _{0}^{i},p_{0}^{i},U_{i},S_{i}),$%
\textit{\ }$i=1,2,$\textit{\ we have} 
\begin{eqnarray}
&&\left\Vert (\phi _{1}-\phi _{2})(t)\right\Vert _{H}^{2}+\left\Vert (\sigma
_{1}-\sigma _{2})(t)\right\Vert _{H}^{2}+\left\Vert
(p_{1}-p_{2})(t)\right\Vert _{H}^{2}  \notag \\
&&+\left\Vert \phi _{1}-\phi _{2}\right\Vert
_{L^{2}(0,T;V_{0})}^{2}+\left\Vert \sigma _{1}-\sigma _{2}\right\Vert
_{L^{2}(0,T;V)}^{2}+\left\Vert p_{1}-p_{2}\right\Vert _{L^{2}(0,T;V)}^{2} 
\notag \\
&\leq &C\left( \left\Vert \phi _{0}^{1}-\phi _{0}^{2}\right\Vert
_{H}^{2}+\left\Vert \sigma _{0}^{1}-\sigma _{0}^{2}\right\Vert
_{H}^{2}+\left\Vert p_{0}^{1}-p_{0}^{2}\right\Vert _{H}^{2}\right.  \notag \\
&&\left. \mbox{ \ \ \ }+\left\Vert U_{1}-U_{2}\right\Vert
_{L^{2}(0,T;H)}^{2}+\left\Vert S_{1}-S_{2}\right\Vert
_{L^{2}(0,T;H)}^{2}\right)  \label{e14}
\end{eqnarray}%
\textit{for all} $t\in \lbrack 0,T].$ \textit{Finally, if }$(\phi
_{0},\sigma _{0},p_{0})\in V_{0}\times V\times V,$\textit{\ the solution has
the supplementary regularity }$(\phi ,\sigma ,p)\in \mathcal{X}_{0}\times 
\mathcal{X}\times \mathcal{X},$ \textit{and satisfies the estimate }%
\begin{eqnarray}
&&\left\Vert \phi \right\Vert^2_{\mathcal{X}_{0}}+\left\Vert \sigma
\right\Vert^2_{\mathcal{X}}+\left\Vert p\right\Vert^2_{\mathcal{X}}
\notag \\
&\leq &C(\left\Vert \phi _{0}\right\Vert _{V_{0}}^{2}+\left\Vert \sigma
_{0}\right\Vert _{V}^{2}+\left\Vert p_{0}\right\Vert _{V}^{2}+\left\Vert
U\right\Vert _{L^{2}(0,T;H)}^{2}+\left\Vert S\right\Vert
_{L^{2}(0,T;H)}^{2}+1).  \label{e13}
\end{eqnarray}

\medskip

The proof of this result can be found in Ref.~\cite{CGLMRR-1}.

\medskip

\noindent \textbf{Theorem 3.3. }\textit{Assume } (\ref{e12-0}), (\ref{e13-0}),%
\textit{\ and }(\ref{5-0}).\textit{\ Then, there exists at least one
solution }$(U^{\ast },S^{\ast })\in \mathcal{U}_{ad}$ \textit{to problem }$%
(P),$ \textit{with the corresponding optimal state} $(\phi ^{\ast },\sigma
^{\ast },p^{\ast })$ \textit{solving} (\ref{1})-(\ref{5}) in the sense of
Definition 3.1.

\medskip

\noindent \textbf{Proof.} It is easy to check that the \pierbis{functional\textbf{\ }$%
J$} is nonnegative, so that it has an infimum $d\geq 0.$ We consider a
minimizing sequence $(U_{n},S_{n})\in \mathcal{U}_{ad}$ such that%
\begin{equation}
d\leq J(U_{n},S_{n})\leq d+\frac{1}{n},\mbox{ for }n\geq 1.  \label{11}
\end{equation}%
The state system corresponding to $(U_{n},S_{n})$ has, according to Theorem
3.2, a unique solution $(\phi _{n},\sigma _{n},p_{n})$ satisfying the
estimates (\ref{e12-1}), (\ref{e12-00}), and (\ref{e13-3}). Selecting a
subsequence, still denoted by $n,$ we can write 
\begin{equation*}
U_{n}\rightarrow U^{\ast },\mbox{ }S_{n}\rightarrow S^{\ast }%
\mbox{
weak-star in }L^{\infty }(Q_{T}),
\end{equation*}%
\begin{equation*}
\phi _{n}\rightarrow \phi ^{\ast }\mbox{ weakly in }W^{1,2}(0,T;V_{0}^{%
\prime })\cap L^{2}(0,T;V_{0}),\mbox{ weak-star in }L^{\infty }(Q_{T}),
\end{equation*}%
\begin{equation*}
\sigma _{n}\rightarrow \sigma ^{\ast }\mbox{ weakly in }W^{1,2}(0,T;V^{%
\prime })\cap L^{2}(0,T;V),\mbox{ weak-star in }L^{\infty }(Q_{T}),
\end{equation*}%
\begin{equation*}
p_{n}\rightarrow p^{\ast }\mbox{ weakly in }W^{1,2}(0,T;V^{\prime })\cap
L^{2}(0,T;V),\mbox{ weak-star in }L^{\infty }(Q_{T}).
\end{equation*}%
Applying Lions' compactness theorem (see Ref.~\cite{Lions-69}, p. 58) we deduce
that 
\begin{equation}
\phi _{n}\rightarrow \phi ^{\ast }\mbox{, }\sigma _{n}\rightarrow \sigma
^{\ast },\mbox{ }p_{n}\rightarrow p^{\ast }\mbox{ strongly in }L^{2}(0,T;H).
\label{Lions}
\end{equation}%
By Ascoli-Arzel\`{a} theorem we also obtain that 
\begin{equation}
\phi _{n}(t)\rightarrow \phi ^{\ast }(t)\mbox{ strongly in }V_{0}^{\prime },%
\mbox{ uniformly in }t\in \lbrack 0,T],  \label{AA1}
\end{equation}%
\begin{equation}
\sigma _{n}(t)\rightarrow \sigma ^{\ast }(t)\mbox{ strongly in }V^{\prime },%
\mbox{ uniformly in }t\in \lbrack 0,T],  \label{AA2}
\end{equation}%
\begin{equation}
p_{n}(t)\rightarrow p^{\ast }(t)\mbox{ strongly in }V^{\prime },%
\mbox{
uniformly in }t\in \lbrack 0,T],  \label{AA3}
\end{equation}%
whence $\phi ^{\ast }(0)=\phi _{0},$ $\sigma ^{\ast }(0)=\sigma _{0},$ $%
p^{\ast }(0)=p_{0}.$ In addition, $\phi ^{\ast },$ $\sigma ^{\ast },$ $%
p^{\ast }\in L^{\infty }(Q_{T})$ and satisfy (\ref{e12-1}), (\ref{e13-2}),
and (\ref{e13-3}).

Using the boundedness and the Lipschitz continuity of $F^{\prime }$ and $%
h^{\prime }$ in the interval $[0,1],$ we derive that 
\begin{equation*}
F^{\prime }(\phi _{n})\rightarrow F^{\prime }(\phi ^{\ast }),\mbox{ }%
h^{\prime }(\phi _{n})\rightarrow h^{\prime }(\phi ^{\ast })%
\mbox{ strongly
in }L^{2}(Q_{T}),\mbox{ weak-star in }L^{\infty }(Q_{T}).
\end{equation*}%
By (\ref{7}) and (\ref{8}) we obtain%
\begin{equation*}
m(\sigma _{n})\rightarrow m(\sigma ^{\ast })\mbox{ strongly in }L^{2}(Q_{T}),%
\mbox{ weak-star in }L^{\infty }(Q_{T}),
\end{equation*}%
whence%
\begin{equation*}
m(\sigma _{n})h^{\prime }(\phi _{n})\rightarrow m(\sigma ^{\ast })h^{\prime
}(\phi ^{\ast })\mbox{ strongly in }L^{2}(Q_{T}),\mbox{ }
\end{equation*}%
\begin{equation*}
\phi _{n}\sigma _{n}\rightarrow \phi ^{\ast }\sigma ^{\ast }%
\mbox{ strongly
in }L^{2}(Q_{T}),\mbox{ }
\end{equation*}%
and 
\begin{gather*}
U_{n}h^{\prime }(\phi _{n})\rightarrow U^{\ast }h^{\prime }(\phi ^{\ast })%
\mbox{ weakly in }L^{2}(Q_{T}), \\
\pg{S_{n}\phi _{n}\rightarrow S^{\ast } \phi ^{\ast } \mbox{ weakly in
}L^{2}(Q_{T}). }
\end{gather*}%
Then, we can pass to the limit in the variational formulations (\ref{e9-1})-(%
\ref{e10}) for the solution $(\phi _{n},\sigma _{n},p_{n})$ and deduce that $%
(\phi ^{\ast },\sigma ^{\ast },p^{\ast })$ is the solution to %
\pier{(\ref{1})-(\ref{5})} corresponding to $(U^{\ast },S^{\ast }).$

From (\ref{e12-00}) we see that $\phi _{n}(T),$ $p_{n}(T)$ are uniformly
bounded in $H,$ so that $\phi _{n}(T)\rightarrow \phi ^{\ast }(T)$ and $%
p_{n}(T)\rightarrow p^{\ast }(T)$ weakly in $H$ by (\ref{AA1}) and (\ref{AA3}%
), and consequently%
\begin{equation}
k_{3}\int_{\Omega }\phi _{n}(T)dx\rightarrow k_{3}\int_{\Omega }\phi ^{\ast
}(T)dx,\mbox{ \ }k_{5}\int_{\Omega }p_{n}(T)dx\rightarrow k_{5}\int_{\Omega
}p^{\ast }(T)dx.  \label{12}
\end{equation}

Also, with the help of (\ref{Lions}) and the Lipschitz continuity of the
positive part, it turns out that 
\begin{equation}
\left( \int_{\Omega }p_{n}(\cdot )dx-p_{\Omega }(\cdot )\right)
^{+}\rightarrow \left( \int_{\Omega }p^{\ast }(\cdot )dx-p_{\Omega }(\cdot
)\right) ^{+}\mbox{ strongly in }L^{2}(0,T).  \label{12-0}
\end{equation}%
We aim at showing that $J(U^{\ast },S^{\ast })\leq \liminf\limits_{n\rightarrow
\infty }J(U_{n},S_{n})=\lim\limits_{n\rightarrow \infty }J(U_{n},S_{n})=d$
by (\ref{11}), whence $J(U^{\ast },S^{\ast })=d$ and so $(U^{\ast },S^{\ast
})$ is an optimal control in $(P).$ Indeed, the second and the last two
terms in $J(U_{n},S_{n})$ are weakly lower semicontinuous and 
\begin{eqnarray*}
&&\frac{k_{2}}{2}\int_{\Omega }(\phi ^{\ast }(T)-\phi _{\Omega })^{2}dx+%
\frac{k_{6}}{2}\int_{Q_{T}}(U^{\ast })^{2}dxdt+\frac{k_{7}}{2}%
\int_{Q_{T}}(S^{\ast })^{2}dxdt \\
&\leq &\liminf\limits_{n\rightarrow \infty }\left( \frac{k_{2}}{2}%
\int_{\Omega }(\phi _{n}(T)-\phi _{\Omega })^{2}dx+\frac{k_{6}}{2}%
\int_{Q_{T}}U_{n}^{2}dxdt+\frac{k_{7}}{2}\int_{Q_{T}}S_{n}^{2}dxdt\right) .
\end{eqnarray*}%
On the other hand, we point out that the other four therms in $%
J(U_{n},S_{n}) $ converge, as previously specified. This concludes the
proof. \hfill $\square $

\section{Optimality conditions}\label{optcon}

\setcounter{equation}{0}

Let $(U^{\ast },S^{\ast })\in U_{ad}$\ be an optimal control in $(P)$\ with
the corresponding state $(\phi ^{\ast },\sigma ^{\ast },p^{\ast })$ and let $%
(u,s)\in L^{2}(0,T;H)\times L^{2}(0,T;H).$ We introduce the \pierbis{linearized} 
system 
\begin{align}
&Y_{t}-\lambda \Delta Y+(F^{\prime \prime }(\phi ^{\ast })-m(\sigma ^{\ast
})h^{\prime \prime }(\phi ^{\ast })+U^{\ast }h^{\prime \prime }(\phi ^{\ast
}))Y-m^{\prime }(\sigma ^{\ast })h^{\prime }(\phi ^{\ast })Z  \notag \\
=\ &-uh^{\prime }(\phi ^{\ast }),\mbox{ in }Q_{T},  \label{13}
\end{align}%
\begin{equation}
Z_{t}-\eta \Delta Z+(\gamma _{h}+\gamma _{ch}\phi ^{\ast })Z + \pg{(\gamma
_{ch}\sigma ^{\ast }+ S^\ast - S_{ch}) Y= - s \phi^\ast} \mbox{ in }Q_{T},
\label{14}
\end{equation}%
\begin{equation}
P_{t}-D\Delta P+\gamma _{p}P-\alpha _{ch}Y=0,\mbox{ in }Q_{T},  \label{15}
\end{equation}%
\begin{equation}
Y=0,\frac{\partial Z}{\partial \nu }=\frac{\partial P}{\partial \nu }=0,%
\mbox{ on }\Sigma _{T},  \label{16}
\end{equation}%
\begin{equation}
Y(0)=0,\mbox{ }Z(0)=0,\mbox{ }P(0)=0,\mbox{ in }\Omega .  \label{17}
\end{equation}

\medskip

\noindent \textbf{Proposition 4.1. }\textit{The system }(\ref{13})-(\ref{17}%
) \textit{has a unique \pier{strong} solution} $(Y,Z,P)\in \mathcal{X}%
_{0}\times \mathcal{X}\times \mathcal{X},$ \textit{satisfying} 
\begin{equation}
\left\Vert Y\right\Vert^2_{\mathcal{X}_{0}}+\left\Vert Z\right\Vert^2_{%
\mathcal{X}}+\left\Vert P\right\Vert^2_{\mathcal{X}}\leq C(\left\Vert
u\right\Vert _{L^{2}(0,T;H)}^{2}+\left\Vert s\right\Vert _{L^{2}(0,T;H)}^{2})%
\pier{,}  \label{17-0}
\end{equation}
\pier{\textit{where the spaces $\mathcal{X}_{0}$ and $\mathcal{X}$ are defined in}
(\ref{10}).}

\medskip

\noindent \textbf{Proof. }We note that equations\textbf{\ }(\ref{13})-(\ref%
{15}) are linear parabolic and have the coefficients of $Y,$ $Z$ and $P$ in $%
L^{\infty }(Q_{T})$ and the source terms in $L^{2}(0,T;H),$ according to (%
\ref{6})-(\ref{8}) and to the property $0\leq \phi ^{\ast }\leq 1$ a.e. in $%
Q_{T}.$ The initial data are zero, then smooth. Therefore, by the general
results concerning the existence and uniqueness of the solutions to linear
parabolic systems (cf. Ref.~\cite{Lions-61}), we deduce that there is a unique
triplet $(Y,Z,P)\in \mathcal{X}_{0}\times \mathcal{X}\times \mathcal{X}$
solving (\ref{13})-(\ref{15})\textbf{\ }and satisfying estimate\textbf{\ }(%
\ref{17-0})\textbf{.}\hfill $\square $

\medskip

Next, let $(U^{\ast },S^{\ast })\in U_{ad}$\ be an optimal control in $(P)$\
with the corresponding state $(\phi ^{\ast },\sigma ^{\ast },p^{\ast })$.
Let $\mu \in (0,1)$ and define%
\begin{equation}
U^{\mu }=U^{\ast }+\mu u,\mbox{ }S^{\mu }=S^{\ast }+\mu s,  \label{19}
\end{equation}%
where 
\begin{equation}
(u,s)=(\bar{u}-U^{\ast },\bar{s}-S^{\ast }),\mbox{ for some }(\bar{u},\bar{s}%
)\in \mathcal{U}_{ad}.  \label{18}
\end{equation}
We set%
\begin{equation}
Y^{\mu }=\frac{\phi ^{\mu }-\phi ^{\ast }}{\mu }-Y,\mbox{ }Z^{\mu }=\frac{%
\sigma ^{\mu }-\sigma ^{\ast }}{\mu }-Z,\mbox{ }P^{\mu }=\frac{p^{\mu
}-p^{\ast }}{\mu }-P,  \label{20}
\end{equation}%
where $(\phi ^{\mu },\sigma ^{\mu },p^{\mu })$ is the solution to (\ref{e9-1}%
)-(\ref{e10-0}) corresponding to $(U^{\mu },S^{\mu })\in \mathcal{U}_{ad}$
and $(Y,Z,P)$ solves (\ref{13})-(\ref{17}) with $(u,s)$ defined by (\ref{18}%
).

\medskip

\noindent \textbf{Theorem 4.2. }\textit{The following convergence properties
hold }%
\begin{equation*}
Y^{\mu }\rightarrow 0\mbox{ \textit{strongly in} }C([0,T];H)\cap
L^{2}(0,T;V_{0}),
\end{equation*}%
\begin{equation*}
Z^{\mu }\rightarrow 0\mbox{ \textit{strongly in} }C([0,T];H)\cap
L^{2}(0,T;V),
\end{equation*}%
\begin{equation*}
P^{\mu }\rightarrow 0\mbox{ \textit{strongly in }}C([0,T];H)\cap L^{2}(0,T;V)
\end{equation*}%
\textit{and so they show that} (\ref{13})-(\ref{17}) \textit{is precisely
the system of first order variations related to }(\ref{1})-(\ref{5}).

\medskip

\noindent \textbf{Proof. }We write the system satisfied by $(Y^{\mu },Z^{\mu
},P^{\mu }):$%
\begin{equation}
Y_{t}^{\mu }-\lambda \Delta Y^{\mu }+a_{1}Y^{\mu }+a_{2}Y+a_{3}Z^{\mu
}+a_{4}Z=\pg{-u(h'(\phi^\mu )-h'(\phi^*))} ,\mbox{ in }Q_{T},  \label{21}
\end{equation}%
\begin{equation}
Z_{t}^{\mu }-\eta \Delta Z^{\mu }+b_{1}Y^{\mu }+b_{2}Y+b_{3}Z^{\mu }= \pg{-
s ( \phi^\mu - \phi^*)} ,\mbox{
in }Q_{T},  \label{22}
\end{equation}%
\begin{equation}
P_{t}^{\mu }-D\Delta P^{\mu }+\gamma _{p}P^{\mu }-\alpha _{ch}Y^{\mu }=0,%
\mbox{ in }Q_{T},  \label{23}
\end{equation}%
\begin{equation}
Y^{\mu }=0,\frac{\partial Z^{\mu }}{\partial \nu }=\frac{\partial P^{\mu }}{%
\partial \nu }=0,\mbox{ on }\Sigma _{T},  \label{24}
\end{equation}%
\begin{equation}
Y^{\mu }(0)=0,\mbox{ }Z^{\mu }(0)=0,\mbox{ }P^{\mu }(0)=0,\mbox{ in }\Omega ,
\label{25}
\end{equation}%
where 
\begin{equation}
a_{1}:=F^{\prime \prime }(\phi _{int}^{\mu })-h^{\prime \prime }(\phi
_{int,\mu })(m(\sigma ^{\ast })-U^{\ast }),  \label{26}
\end{equation}%
\begin{equation}
a_{2}:=F^{\prime \prime }(\phi _{int}^{\mu })-F^{\prime \prime }(\phi ^{\ast
})-(h^{\prime \prime }(\phi _{int,\mu })-h^{\prime \prime }(\phi ^{\ast
}))(m(\sigma ^{\ast })-U^{\ast }),  \label{27}
\end{equation}%
\begin{equation}
a_{3}:=-m^{\prime }(\sigma _{int}^{\mu })h^{\prime }(\phi ^{\mu }),
\label{28}
\end{equation}%
\begin{equation}
a_{4}:=m^{\prime }(\sigma ^{\ast })h^{\prime }(\phi ^{\ast })-m^{\prime
}(\sigma _{int}^{\mu })h^{\prime }(\phi ^{\mu }),  \label{29}
\end{equation}%
\begin{equation}
b_{1}:= \pg{\gamma _{ch}\sigma ^{\mu}+ S^\ast - S_{ch}} ,\mbox{ }%
b_{2}:=\gamma _{ch}(\sigma ^{\mu }-\sigma ^{\ast }),\mbox{ }b_{3}:=\gamma
_{h}+\gamma _{ch}\phi ^{\ast }.  \label{30}
\end{equation}%
Here, $\phi _{int}^{\mu },$ $\sigma _{int}^{\mu }$ and $\phi _{int,\mu }$
are \pier{measurable functions (see Appendix A)} occurring in the Taylor
expansions of $F^{\prime }(\phi ^{\mu }),$ $m(\sigma ^{\mu }),$ $h^{\prime
}(\phi ^{\mu }),$ 
\pier{indeed, having
the form 
\[
\phi _{int}^{\mu }=\theta _{\mu }\phi ^{\ast }+(1-\theta _{\mu })\phi ^{\mu
},\quad \pier{\mbox{for some function $\theta _{\mu }$ taking values in } (0,1),}
\]with similar expressions for $\sigma _{int}^{\mu }$ and $\phi _{int,\mu }$.}
Since $U^{\mu }\rightarrow U^{\ast }$ and $S^{\mu }\rightarrow S^{\ast }$
strongly in $L^{2}(Q_{T})$ as $\mu \rightarrow 0,$ then by estimate (\ref%
{e14}) we see that 
\begin{equation}
\phi ^{\mu }\rightarrow \phi ^{\ast },\mbox{ }\sigma ^{\mu }\rightarrow
\sigma ^{\ast },\mbox{ }p^{\mu }\rightarrow p^{\ast }\mbox{ strongly in }%
C([0,T];H),  \label{30-0}
\end{equation}%
whence 
\begin{equation}
\phi _{int}^{\mu }\rightarrow \phi ^{\ast },\mbox{ }\sigma _{int}^{\mu
}\rightarrow \sigma ^{\ast },\mbox{ }\phi _{int,\mu }\rightarrow \phi ^{\ast
}\mbox{ strongly in }C([0,T];H).  \label{30-1}
\end{equation}%
We test (\ref{21}) by $Y^{\mu },$ (\ref{22}) by $Z^{\mu }$, (\ref{23}) by $%
P^{\mu }$ and sum up all of them, obtaining%
\begin{eqnarray}
&&\frac{1}{2}\left\Vert Y^{\mu }(t)\right\Vert _{H}^{2}+\lambda
\int_{0}^{t}\left\Vert \nabla Y^{\mu }(\tau )\right\Vert _{H}^{2}d\tau +%
\frac{1}{2}\left\Vert Z^{\mu }(t)\right\Vert _{H}^{2}+\eta
\int_{0}^{t}\left\Vert \nabla Z^{\mu }(\tau )\right\Vert _{H}^{2}d\tau 
\mbox{
\ \ \ \ } \qquad  \notag \\
&&+\frac{1}{2}\left\Vert P^{\mu }(t)\right\Vert
_{H}^{2}+D\int_{0}^{t}\left\Vert \nabla P^{\mu }(\tau )\right\Vert
_{H}^{2}d\tau +\gamma _{p}\int_{0}^{t}\left\Vert P^{\mu }(\tau )\right\Vert
_{H}^{2}d\tau  \notag \\
&\leq &\int_{Q_{t}}|a_{1}||Y^{\mu }|^{2}dxd\tau
+\int_{Q_{t}}|a_{2}||Y||Y^{\mu }|dxd\tau +\int_{Q_{t}}|a_{3}||Z^{\mu
}||Y^{\mu }|dxd\tau  \notag \\
&&+\int_{Q_{t}}|a_{4}||Z||Y^{\mu }|dxd\tau +\int_{Q_{t}}|b_{1}||Y^{\mu
}||Z^{\mu }|dxd\tau +\int_{Q_{t}}|b_{2}||Y||Z^{\mu }|dxd\tau  \notag \\
&&+\int_{Q_{t}}|b_{3}||Z^{\mu }|^{2}dxd\tau +\int_{Q_{t}}|\alpha
_{ch}||Y^{\mu }||P^{\mu }|dxd\tau  \notag \\
&&+\pg{{}\int_{Q_{t}}|u||h'(\phi^{\mu })-h'( \phi^*)||Y^{\mu }|dxd\tau
+\int_{Q_{t}}|s||\phi^{\mu }- \phi^*||Z^{\mu }|dxd\tau} .  \label{31}
\end{eqnarray}%
We see that $a_{1},$ $a_{3},$ $b_{1},$ $b_{3}$ are uniformly bounded in $%
L^{\infty }(Q_{T})$ \pierbis{and, by \eqref{18}, we have that $|u| \le U_{\max}$, 
$|s| \le S_{\max}$ a.e in $Q_T$.}
Then, by the Young inequality, it follows that 
\begin{eqnarray*}
&&\int_{Q_{t}}|a_{1}||Y^{\mu }|^{2}dxd\tau +\int_{Q_{t}}|a_{3}||Z^{\mu
}||Y^{\mu }|dxd\tau +\int_{Q_{t}}|b_{1}||Y^{\mu }||Z^{\mu }|dxd\tau \\
&&+\int_{Q_{t}}|b_{3}||Z^{\mu }|^{2}dxd\tau +\int_{Q_{t}}|\alpha
_{ch}||Y^{\mu }||P^{\mu }|dxd\tau \\
&\leq &C\left( \int_{Q_{t}}|Y^{\mu }|^{2}dxd\tau +\int_{Q_{t}}|Z^{\mu
}|^{2}dxd\tau +\int_{Q_{t}}|P^{\mu }|^{2}dxd\tau \right) .
\end{eqnarray*}%
Next, by applying the H\"{o}lder inequality and the continuous embedding $%
V_{0}\subset L^{4}(\Omega ),$ we easily infer that 
\begin{eqnarray*}
&&\int_{Q_{t}}|a_{2}||Y||Y^{\mu }|dxd\tau \leq \int_{0}^{t}\left\Vert
a_{2}(\tau )\right\Vert _{H}\left\Vert Y(\tau )\right\Vert _{4}\left\Vert
Y^{\mu }(\tau )\right\Vert _{4}d\tau \\
&\leq &C\int_{0}^{t}\left\Vert a_{2}(\tau )\right\Vert _{H}\left\Vert Y(\tau
)\right\Vert _{V_{0}}\left\Vert Y^{\mu }(\tau )\right\Vert _{V_{0}}d\tau \\
&\leq &\frac{\lambda }{4}\int_{0}^{t}\left\Vert \nabla Y^{\mu }(\tau
)\right\Vert _{H}^{2}d\tau +C\left\Vert Y\right\Vert _{L^{\infty
}(0,T;V_{0})}^{2}\int_{0}^{t}\left\Vert a_{2}(\tau )\right\Vert _{H}^{2}d\tau
\\
&\leq &\frac{\lambda }{4}\int_{0}^{t}\left\Vert \nabla Y^{\mu }(\tau
)\right\Vert _{H}^{2}d\tau +C\int_{0}^{t}\left\Vert a_{2}(\tau )\right\Vert
_{H}^{2}d\tau ,
\end{eqnarray*}%
where we used (\ref{17-0}). The \pg{following} remaining terms on the
right-hand side of (\ref{31}) can be treated in the same way, in order to
find out that 
\begin{eqnarray*}
&&\int_{Q_{t}}|a_{4}||Z||Y^{\mu }|dxd\tau +\int_{Q_{t}}|b_{2}||Y||Z^{\mu
}|dxd\tau \\
&\leq &\frac{\lambda }{4}\int_{0}^{t}\left\Vert \nabla Y^{\mu }(\tau
)\right\Vert _{H}^{2}d\tau +C\int_{0}^{t}\left\Vert a_{4}(\tau )\right\Vert
_{H}^{2}d\tau +\frac{\eta }{2}\int_{0}^{t}\left\Vert \nabla Z^{\mu }(\tau
)\right\Vert _{H}^{2}d\tau \\
&&\pier{{}+\frac{\eta }{2}\int_{0}^{t}\left\Vert Z^{\mu }(\tau )\right\Vert
_{H}^{2}d\tau} +C\int_{0}^{t}\left\Vert b_{2}(\tau )\right\Vert
_{H}^{2}d\tau .
\end{eqnarray*}%
Going back to (\ref{31}), we obtain 
\begin{align}
&\frac{1}{2}\left\Vert Y^{\mu }(t)\right\Vert _{H}^{2}+\frac{\lambda }{2}%
\int_{0}^{t}\left\Vert \nabla Y^{\mu }(\tau )\right\Vert _{H}^{2}d\tau +%
\frac{1}{2}\left\Vert Z^{\mu }(t)\right\Vert _{H}^{2}+\frac{\eta }{2}%
\int_{0}^{t}\left\Vert \nabla Z^{\mu }(\tau )\right\Vert _{H}^{2}d\tau \notag \\
&{}+\frac{1}{2}\left\Vert P^{\mu }(t)\right\Vert
_{H}^{2}+D\int_{0}^{t}\left\Vert \nabla P^{\mu }(\tau )\right\Vert
_{H}^{2}d\tau +\gamma _{p}\int_{0}^{t}\left\Vert P^{\mu }(\tau )\right\Vert
_{H}^{2}d\tau  \notag \\
\leq\ &C\left( \int_{0}^{t}\left\Vert Y^{\mu }(\tau )\right\Vert
_{H}^{2}d\tau +\int_{0}^{t}\left\Vert Z^{\mu }(\tau )\right\Vert
_{H}^{2}d\tau +\int_{0}^{t}\left\Vert P^{\mu }(\tau )\right\Vert
_{H}^{2}d\tau \right)  \notag \\
&{}+C\left( \int_{0}^{t}\left\Vert a_{2}(\tau )\right\Vert _{H}^{2}d\tau
+\int_{0}^{t}\left\Vert a_{4}(\tau )\right\Vert _{H}^{2}d\tau
+\int_{0}^{t}\left\Vert b_{2}(\tau )\right\Vert _{H}^{2}d\tau %
 \pg{{}}\right)   \notag \\
&\pierbis{{} +C \int_{0}^{t}\left\Vert (\phi^\mu - \phi^*)(\tau )\right\Vert
_{H}^{2}d\tau}, \label{32}
\end{align}%
\pg{by using the Lipschitz continuity of the function $h'$.} An application
of the Gronwall lemma leads to%
\begin{eqnarray}
&&\left\Vert Y^{\mu }(t)\right\Vert _{H}^{2}+\int_{0}^{t}\left\Vert \nabla
Y^{\mu }(\tau )\right\Vert _{H}^{2}d\tau +\left\Vert Z^{\mu }(t)\right\Vert
_{H}^{2}+\int_{0}^{t}\left\Vert \nabla Z^{\mu }(\tau )\right\Vert
_{H}^{2}d\tau  \notag \\
&&+\left\Vert P^{\mu }(t)\right\Vert _{H}^{2}+\int_{0}^{t}\left\Vert P^{\mu
}(\tau )\right\Vert _{V}^{2}d\tau  \notag \\
&\leq &C\left( \left\Vert a_{2}\right\Vert _{L^{2}(Q_{T})}^{2}+\left\Vert
a_{4}\right\Vert _{L^{2}(Q_{T})}^{2}+\left\Vert b_{2}\right\Vert
_{L^{2}(Q_{T})}^{2} \right. \notag \\ && \left. +\left\Vert \phi^\mu - \phi^*\right\Vert
_{L^{2}(Q_{T})}^{2} \right) ,\mbox{ for all }t\in \lbrack 0,T].  \label{33}
\end{eqnarray}%
Recalling (\ref{30-0}), (\ref{30-1}) and the boundedness and Lipschitz
continuity of $F^{\prime \prime },$ $h^{\prime },$ $h^{\prime \prime }$ in $%
[0,1]$ and of $m^{\prime }$ in $\mathbb{R},$ it is not difficult to check
that $a_{2},$ $a_{4}$, \pg{$b_{2}$ and the last term in \eqref{33}} tend to
zero strongly in $L^{2}(Q_{T}),$ as $\mu \rightarrow 0.$ Thanks to (\ref{33}%
) this implies that 
\begin{equation*}
Y^{\mu }\rightarrow 0,\mbox{ }Z^{\mu }\rightarrow 0,\mbox{ }P^{\mu
}\rightarrow 0\mbox{ strongly in }C([0,T];H)\cap L^{2}(0,T;V).
\end{equation*}%
Therefore, we see that 
\begin{equation*}
\lim_{\mu \rightarrow 0}\frac{\phi ^{\mu }-\phi ^{\ast }}{\mu }=Y,\mbox{ }%
\lim_{\mu \rightarrow 0}\frac{\sigma ^{\mu }-\sigma ^{\ast }}{\mu }=Z,%
\mbox{ 
}\lim_{\mu \rightarrow 0}\frac{p^{\mu }-p^{\ast }}{\mu }=P,
\end{equation*}%
where the limit is understood in the sense of the strong convergence in $%
C([0,T];H)\cap L^{2}(0,T;V).$ This concludes the proof.\hfill $\square $

\medskip

At this point, we can introduce the adjoint system in terms of the adjoint
variables $w,$ $z,$ $q:$%
\begin{align}
&-w_{t}-\lambda \Delta w+(F^{\prime \prime }(\phi ^{\ast })-m(\sigma ^{\ast
})h^{\prime \prime }(\phi ^{\ast })+U^{\ast }h^{\prime \prime }(\phi ^{\ast
}))w  \notag \\
&{}+\pg{(\gamma _{ch}\sigma ^{\ast } + (S^* - S_{ch}) )} z-\alpha _{ch}q
=k_{1}(\phi ^{\ast }-\phi _{Q}),\mbox{ in }Q_{T},  \label{34}
\end{align}%
\begin{equation}
-z_{t}-\eta \Delta z+(\gamma _{h}+\gamma _{ch}\phi ^{\ast })z-m^{\prime
}(\sigma ^{\ast })h^{\prime }(\phi ^{\ast })w=0,\mbox{ in }Q_{T},  \label{35}
\end{equation}%
\begin{equation}
-q_{t}-D\Delta q+\gamma _{p}q=k_{4}\left( \int_{\Omega }p^{\ast }(\cdot
)dx-p_{\Omega }(\cdot )\right) ^{+},\mbox{ in }Q_{T},  \label{36}
\end{equation}%
\begin{equation}
w=0,\frac{\partial z}{\partial \nu }=\frac{\partial q}{\partial \nu }=0,%
\mbox{ on }\Sigma _{T},  \label{37}
\end{equation}%
\begin{equation}
w(T)=k_{2}(\phi ^{\ast }(T)-\phi _{\Omega })+k_{3},\mbox{ }z(T)=0,\mbox{ }%
q(T)=k_{5},\mbox{ in }\Omega ,  \label{38}
\end{equation}%
where $k_{i},$ $i=1,...,5$, \pierbis{and} $\phi _{Q},$ $\phi _{\Omega },$ $p_{\Omega }$
are the coefficients and the functions involved in the cost functional (\ref%
{J}).

\medskip

\noindent \textbf{Proposition 4.3. }\textit{The system }(\ref{34})-(\ref{38}%
) \textit{has a unique variational solution} $(w,z,q)\in X_{0}\times X\times
X,$ \textit{satisfying} 
\begin{equation}
\left\Vert w\right\Vert _{X_{0}}+\left\Vert z\right\Vert _{X}+\left\Vert
q\right\Vert _{X}\leq C,  \label{39}
\end{equation}%
\textit{where }$C$\textit{\ is a constant depending on the norms} $%
\left\Vert F^{\prime \prime }\right\Vert _{L^{\infty }([0,1])},$ $\left\Vert
m\right\Vert _{L^{\infty }(\mathbb{R})},$ $\left\Vert h^{\prime \prime
}\right\Vert _{L^{\infty }([0,1])},$ $\left\Vert \sigma ^{\ast
}\right\Vert _{\infty },$ $\left\Vert \phi ^{\ast }(T)-\phi _{\Omega
}\right\Vert _{L^{2}(\Omega )},$  $\left\Vert \phi ^{\ast }-\phi _{Q}\right\Vert
_{L^{2}(Q_{T})},$ $\left\Vert m^{\prime }\right\Vert _{L^{\infty }(\mathbb{R}%
)},$ $\left\Vert h^{\prime }\right\Vert _{L^{\infty }([0,1])},$ $\left\Vert
\left( \int_{\Omega }p^{\ast }(\cdot )dx-p_{\Omega }(\cdot )\right)
^{+}\right\Vert _{L^{2}(0,T)}$, $U_{\max}$, $S_{\max}$,  $\Omega ,$ $T,$ \textit{the problem
parameters and all} $k_{i},$ $i=1,...,5.$

\medskip

\noindent \textbf{Proof. }Using the transformation $t\mapsto T-t,$ we can
rewrite (\ref{34})-(\ref{38}) as a linear parabolic system with initial
conditions in $H$ and with the coefficients of $w,$ $z$ and $q$ in $%
L^{\infty }(Q_{T}).$ Moreover, the source terms are in $L^{2}(0,T;H),$
according to (\ref{6})-(\ref{8}) and $0\leq \phi ^{\ast }\leq 1$ a.e. in $%
Q_{T}.$ Therefore, by the general results concerning the existence and
uniqueness of the solutions to linear parabolic systems (cf. Ref.~\cite{Lions-61}%
), we deduce that there is a unique triplet $(w,z,q)\in X_{0}\times X\times
X $ solving (\ref{34})-(\ref{38})\textbf{\ }and satisfying estimate\textbf{\ 
}(\ref{39})\textbf{.}\hfill $\square $

\medskip

We introduce the convex closed subsets of $L^{2}(Q_{T})$: 
\begin{eqnarray*}
K_{1} &:&=\{u\in L^{2}(Q_{T});\mbox{ }0\leq u\leq U_{\max }\mbox{ a.e. in }%
Q_{T}\},\mbox{ } \\
K_{2} &:&=\{s\in L^{2}(Q_{T});\mbox{ }0\leq s\leq S_{\max }\mbox{ a.e. in }%
Q_{T}\}.
\end{eqnarray*}%
We also specify that the normal cone at $U^{\ast }$ to $K_{1}$ is given by%
\begin{equation*}
N_{K_{1}}(U^{\ast })=\{\xi \in L^{2}(Q_{T});\mbox{ }(\xi ,U^{\ast }-\bar{u}%
)\geq 0\mbox{ for all }\bar{u}\in K_{1}\}.
\end{equation*}%
Similarly, we denote by $N_{K_{2}}(S^{\ast })$ the normal cone at $S^{\ast }$
to $K_{2}.$

\medskip

\noindent \textbf{Theorem 4.4. }\textit{Let }$(U^{\ast },S^{\ast })\in
U_{ad} $\textit{\ be an optimal control for }$(P)$\textit{\ with the
corresponding state }$(\phi ^{\ast },\sigma ^{\ast },p^{\ast }).$ \textit{%
Then, the first order optimality conditions}%
\begin{equation}
(k_{6}I+N_{K_{1}})U^{\ast }\ni h^{\prime }(\phi ^{\ast })w,  \label{40}
\end{equation}%
\begin{equation}
(k_{7}I+N_{K_{2}})S^{\ast }\ni \phi ^{\ast }z  \label{41}
\end{equation}%
\textit{hold true, }$w$\textit{\ and }$z$\textit{\ being the first and third
components of the solution }$(w,z,q)\in X_{0}\times X\times X$\textit{\ to
the system} (\ref{34})-(\ref{38}).

\medskip

\noindent \textbf{Proof. }Since $(U^{\ast },S^{\ast })\in U_{ad}$\ is an
optimal control for $(P),$ then we have that 
\begin{equation}
J(U^{\mu },S^{\mu })\geq J(U^{\ast },S^{\ast }),  \label{50}
\end{equation}%
where $U^{\mu }$ and $S^{\mu }$ were defined in (\ref{19}), namely, $U^{\mu
}=U^{\ast }+\mu u,$ $S^{\mu }=S^{\ast }+\mu s,$ with $(u,s)=(\bar{u}-U^{\ast
},\bar{s}-S^{\ast })$ and $(\bar{u},\bar{s})$ arbitrary in $\mathcal{U}%
_{ad}. $ Replacing the expressions of $U^{\mu }$ and $S^{\mu }$ in (\ref{50}%
), by a straightforward calculation we obtain the optimality conditions in
terms of the solution $(Y,Z,P)$ of the system in variations (\ref{13})-(\ref%
{17}) corresponding to $(u,s).$ We\ easily deduce that 
\begin{eqnarray}
&&k_{1}\int_{Q_{T}}(\phi ^{\ast }-\phi _{Q})Ydxdt+k_{2}\int_{\Omega }(\phi
^{\ast }(T)-\phi _{\Omega })Y(T)dx+k_{3}\int_{\Omega }Y(T)dx  \notag \\
&&+k_{4}\int_{0}^{T}\left[ \left( \int_{\Omega }p^{\ast }(t)dx-p_{\Omega
}(t)\right) ^{+}\int_{\Omega }P(t)dx\right] dt+k_{5}\int_{\Omega }P(T)dx 
\notag \\
&&+k_{6}\int_{Q_{T}}U^{\ast }udxdt+k_{7}\int_{Q_{T}}S^{\ast }sdxdt\geq 0. 
\label{51}
\end{eqnarray}%
Now we test (\ref{13}) by $w,$ (\ref{14}) by $z,$ (\ref{15}) by $q$, sum up
and integrate by parts taking into account the  boundary conditions (\ref%
{16}) and (\ref{37}). Using also the initial conditions (\ref{17}) of the
system in variations and the final conditions (\ref{38}) of the adjoint
system, we obtain the equation%
\begin{eqnarray}
&&\int_{\Omega }\left( k_{2}(\phi ^{\ast }(T)-\phi _{\Omega })+k_{3}\right)
Y(T)dx+\int_{\Omega }k_{5}P(T)dx  \notag \\
&&+\int_{0}^{T}\left\langle (-w_{t}-\lambda \Delta w+(F^{\prime \prime }(\phi
^{\ast })-m(\sigma ^{\ast })h^{\prime \prime }(\phi ^{\ast })+U^{\ast
}h^{\prime \prime }(\phi ^{\ast })w)(t),Y(t)\right\rangle _{V_{0}^{\prime
},V_{0}} dt  \notag \\
&&+\int_{0}^{T}\left\langle (\gamma _{ch}\sigma ^{\ast }z-\alpha
_{ch}q)(t),Y(t)\right\rangle _{V_{0}^{\prime },V_{0}}dt  \notag \\
&&+\int_{0}^{T}\left\langle (-z_{t}-\eta \Delta z+(\gamma _{h}+\gamma
_{ch}\phi ^{\ast })z-m^{\prime }(\sigma ^{\ast })h^{\prime }(\phi ^{\ast
})w)(t),Z(t)\right\rangle _{V^{\prime },V}dt  \notag \\
&&+\int_{0}^{T}\left\langle (-q_{t}-D\Delta q+\gamma _{p}q)(t),P(t)\right\rangle
_{V^{\prime },V}dt  \notag \\
&=&\int_{Q_{T}}\left( -uh^{\prime }(\phi ^{\ast })w \pg{{}- s \phi^* z}%
\right) dxdt.  \label{52}
\end{eqnarray}%
Using again the adjoint system (\ref{34})-(\ref{36}) in (\ref{52}) we obtain 
\begin{eqnarray}
&&k_{1}\int_{Q_{T}}(\phi ^{\ast }-\phi _{Q})Ydxdt+k_{2}\int_{\Omega }(\phi
^{\ast }(T)-\phi _{\Omega })Y(T)dx+k_{3}\int_{\Omega }Y(T)dx\mbox{ \ \ }
\notag \\
&&+k_{4}\int_{Q_{T}}\left( \int_{\Omega }p^{\ast }(t)dx-p_{\Omega
}(t)\right) ^{+}Pdxdt+k_{5}\int_{\Omega }P(T)dx  \notag \\
&=&\int_{Q_{T}}\left( -uh^{\prime }(\phi ^{\ast })w \pg{{}- s \phi^* z}%
\right) dxdt. \label{53}
\end{eqnarray}%
By comparison between (\ref{51}) and (\ref{53}) we easily deduce%
\begin{equation*}
\int_{Q_{T}}\left( -uh^{\prime }(\phi ^{\ast })w \pg{{}- s \phi^* z}\right)
dxdt+k_{6}\int_{Q_{T}}U^{\ast }udxdt+k_{7}\int_{Q_{T}}S^{\ast }sdxdt\geq 0,
\end{equation*}%
where $(u,s)=(\bar{u}-U^{\ast },\bar{s}-S^{\ast })$ and $(\bar{u},\bar{s})$
is arbitrary in $\mathcal{U}_{ad}.$ Therefore, we can write%
\begin{equation}
\int_{Q_{T}}\left\{ (h^{\prime }(\phi ^{\ast })w-k_{6}U^{\ast })(U^{\ast }-%
\bar{u})+(-k_{7}S^{\ast } \pg{{}+ \phi^* z})(S^{\ast }-\bar{s})\right\}
dxdt\geq 0,  \label{54}
\end{equation}%
for all $(\bar{u},\bar{s})\in \mathcal{U}_{ad}.$ Then, taking $\bar{s}%
=S^{\ast }$ and $\bar{u}$ arbitrary in (\ref{54}) we obtain (\ref{40}).
Finally, setting $\bar{u}=U^{\ast }$ and $\bar{s}$ arbitrary we derive (\ref%
{41}), as claimed. \hfill $\square $

\medskip

\noindent \textbf{Remark 4.5. }We note that the optimality conditions do not
depend on the component $q$ of the solution to the adjoint system. This is
essentially due to the biological meaning of the model, because $p,$ the
third component of the state system, does not influence the evolution of the
other components $\phi $ and $\sigma ,$ but is a result of their
evolution.

\medskip

\noindent \textbf{Remark 4.6. }Let us denote by Proj$_{K}(v)$ the projection
of $v\in L^{2}(Q_{T})$ into 
\begin{equation*}
K=\{\zeta \in L^{2}(Q_{T});\mbox{ }a\leq \zeta \leq b\mbox{ a.e. in }Q_{T}\}
\end{equation*}%
defined by 
\begin{equation}
\mbox{Proj}_{K}(v)=\left\{ 
\begin{array}{l}
a,\mbox{ if }v<a, \\ 
v,\mbox{ if }a\leq v\leq b\mbox{,} \\ 
b,\mbox{ if }v>b.%
\end{array}%
\right.   \label{42}
\end{equation}%
We also recall that if $(I+\kappa N_{K})v^{\ast }\ni v,$ then 
\begin{equation}
v^{\ast }=(I+\kappa N_{K})^{-1}v=\mbox{Proj}_{K}v\mbox{ for all }\kappa >0.
\label{43}
\end{equation}%
Hence, if the coefficients $k_{6}$ and $k_{7}$ are positive, then (\ref{40})
and (\ref{41}) entail that%
\begin{equation}
U^{\ast }=\mbox{Proj}_{K_{1}}\left( \frac{1}{k_{6}}h^{\prime }(\phi ^{\ast
})w\right) ,\mbox{ }S^{\ast }=\mbox{Proj}_{K_{2}}\left( \frac{1}{k_{7}}\,%
\pg{\phi^*z }\right) .  \label{44}
\end{equation}%
The relations (\ref{44}) can be written as 
\begin{equation}
U^{\ast }=\left\{ 
\begin{array}{l}
0,\mbox{ \ \ \ \ \ \ \ \ \ \ \ \ \ on the set }\{h^{\prime }(\phi ^{\ast
})w<0\}, \\[0.2cm]
\frac{1}{k_{6}}h^{\prime }(\phi ^{\ast })w,\mbox{ \ \ on the set \thinspace }%
\{0\leq h^{\prime }(\phi ^{\ast })w\leq k_{6}U_{\max }\}, \\[0.2cm] 
U_{\max },\mbox{ \ \ \ \ \ \ \ \ \ on the set }\{h^{\prime }(\phi ^{\ast
})w>k_{6}U_{\max }\}%
\end{array}%
\right.   \label{45}
\end{equation}%
and 
\begin{equation}
S^{\ast }=\left\{ 
\begin{array}{l}
0,\mbox{ \ \ \ \ \ \ on the set }\{\pg{\phi^*z }<0\}, \\[0.2cm] 
\frac{1}{k_{7}}\,\pg{\phi^* z},\mbox{ \ on the set \thinspace }\{0\leq %
\pg{\phi^*z }\leq k_{7}S_{\max }\}, \\[0.2cm] 
S_{\max },\mbox{ \ \ on the set }\{\pg{\phi^*z }>k_{7}S_{\max }\}.%
\end{array}%
\right.   \label{46}
\end{equation}%
On the other hand, if a coefficient is zero, let us say, if $k_{6}=0$, (\ref%
{40}) becomes%
\begin{equation}
N_{K_{1}}(U^{\ast })\ni h^{\prime }(\phi ^{\ast })w,  \label{47}
\end{equation}%
which implies that 
\begin{equation}
U^{\ast }\left\{ 
\begin{array}{l}
=0,\mbox{ \ \ \ \ \ \ \ \ \ on the set }\{h^{\prime }(\phi ^{\ast })w<0\},
\\[0.2cm] 
\in \lbrack 0,U_{\max }],\mbox{ on the set \thinspace }\{h^{\prime }(\phi
^{\ast })w=0\}, \\[0.2cm] 
=U_{\max },\mbox{ \ \ \ \ on the set }\{h^{\prime }(\phi ^{\ast })w>0\}.%
\end{array}%
\right.   \label{48}
\end{equation}%
A similar result can be deduced for $S^{\ast }$ if $k_{7}=0.$

\medskip

\noindent \textbf{Remark 4.7.} We discuss here how the previous results
change in the case that the controllers $U$ and $S$ are considered time
dependent only. This formulation might be more realistic for the therapy and, hence we will use it to explore the behavior of the optimal control problem in the simulations presented in Section~\ref{simstudy}.
First, the functional $J$ becomes%
\begin{eqnarray*}
J(U,S) &=&\frac{k_{1}}{2}\int_{Q_{T}}(\phi (t,x)-\phi _{Q})^{2}dxdt+\frac{%
k_{2}}{2}\int_{\Omega }(\phi (T,x)-\phi _{\Omega })^{2}dx+k_{3}\int_{\Omega
}\phi (T,x)dx \\
&&+\frac{k_{4}}{2}\int_{0}^{T}\left[ \left( \int_{\Omega }p(t,x)dx-p_{\Omega
}(t)\right) ^{+}\right] ^{2}dt+k_{5}\int_{\Omega }p(T,x)dx \\
&&+\frac{k_{6}}{2}\int_{0}^{T}U^{2}(t)dt+\frac{k_{7}}{2}%
\int_{0}^{T}S^{2}(t)dt
\end{eqnarray*}%
and the minimization problem is
\begin{equation*}
\mbox{Minimize }\{J(U,S);\mbox{ }(U,S)\in \mathcal{U}_{ad}\},
\end{equation*}%
subject to (\ref{1})-(\ref{5}), where 
\begin{equation*}
\mathcal{U}_{ad}=\{(u,s)\in L^{\infty }(0,T)\times L^{\infty }(0,T);\mbox{ }%
0\leq u\leq U_{\max },\mbox{ }0\leq s\leq S_{\max },\mbox{ a.e. in }(0,T)\}.
\end{equation*}%
All existence results in Theorems 3.2, 3.3, 4.2 and Propositions 4.1, 4.3
are conserved. For computing the optimality conditions we consider the sets%
\begin{eqnarray*}
K_{1} &:&=\{u\in L^{2}(0,T);\mbox{ }0\leq u\leq U_{\max }\mbox{ a.e. in }%
(0,T)\},\mbox{ } \\
K_{2} &:&=\{s\in L^{2}(0,T);\mbox{ }0\leq s\leq S_{\max }\mbox{ a.e. in }%
(0,T)\}
\end{eqnarray*}%
and by remaking the calculation from (\ref{53}) we get 
\begin{align}
\int_{0}^{T} & \left\lbrace \left( \int_{\Omega }h^{\prime }(\phi ^{\ast
})wdx-k_{6}U^{\ast }\right) (U^{\ast }-\bar{u}) \right. \notag \\  & \left. +\left( \int_\Omega\phi
^{\ast }zdx-k_{7}S^{\ast }\right) (S^{\ast }-\bar{s})\right\rbrace dt\geq 0. \label{40-0} 
\end{align}%
Thus, (\ref{40})-(\ref{41}) in Theorem 4.4 are replaced by%
\begin{align}
(k_{6}I+N_{K_{1}})U^{\ast } &\ni \int_{\Omega }h^{\prime }(\phi ^{\ast })wdx, \label{40-1} \\
(k_{7}I+N_{K_{2}})S^{\ast } &\ni \int_{\Omega }\phi ^{\ast }zdx.  \label{41-1}
\end{align}

\medskip

\noindent \textbf{Remark 4.8. }For a later use we also recall that from a
condition of optimality written as 
\begin{equation*}
\int_0^T(d_{U}u+d_{S}s)dt=\int_0^T\left\{ d_{U}(\bar{u}-U^{\ast
})+d_{S}(\bar{s}-S^{\ast })\right\} dt\geq 0
\end{equation*}%
we can formally read that%
\begin{equation}
d_{U}=\nabla_{U} J(U^{\ast },S^{\ast }),\mbox{ }d_{S}=\nabla _{S}J(U^{\ast
},S^{\ast }),  \label{54-0}
\end{equation}%
where $\nabla_{U} J (U^{\ast },S^{\ast })$ and $\nabla _{S}J(U^{\ast
},S^{\ast })$ are the derivatives of $J$ with respect to the first and
second variables, respectively, calculated at $(U^{\ast },S^{\ast }).$
Consequently, we have%
\begin{align}
d_{U} &= \nabla _{U}J(U^{\ast },S^{\ast })=k_{6}U^{\ast }-\int_{\Omega
}h^{\prime }(\phi ^{\ast })wdx,  \label{54-1} \\
d_{S} &= \nabla _{S}J(U^{\ast },S^{\ast })=k_{7}S^{\ast }-\int_{\Omega }\phi
^{\ast }zdx.  \label{54-2}
\end{align}

\section{Numerical method}\label{numerics}

\subsection{Spatial discretization}\label{space_disc}

We use Isogeometric Analysis \cite{hughes2005isogeometric}, a recent generalization of Finite Element Analysis, to discretize in space the forward and adjoint problems. Our spatial discretization of the forward problem is based on the following weak form of Eqns.~\eqref{1}--\eqref{3}: find $\phi\in V_0$, $\sigma\in V$, and $p\in V$, such that
\begin{align} 
\label{bwf1} & B_1^\chi(\chi_1,\phi,\sigma,p)=0\quad\text{for all}\quad \chi_1\in V_0, \\
\label{bwf2} & B_2^\chi(\chi_2,\phi,\sigma,p)=0\quad\text{for all}\quad \chi_2\in V,   \\
\label{bwf3} & B_3^\chi(\chi_3,\phi,\sigma,p)=0\quad\text{for all}\quad \chi_3\in V, 
\end{align}
where 
\begin{align} 
\label{wf1} B_1^\chi(\chi_1,\phi,\sigma,p) & = \int_{\Omega} \chi_1 \left[ \phi_t + F^\prime(\phi) + h^\prime(\phi)(U-m(\sigma)) \right]dx + \int_{\Omega} \lambda\nabla \chi_1\cdot  \nabla \phi dx, \\
\label{wf2} B_2^\chi(\chi_2,\phi,\sigma,p) & = \int_{\Omega} \chi_2 \left[ \sigma_t + \gamma_h\sigma + \gamma_{ch}\sigma \phi  - S_h(1-\phi) - (S_c-S)\phi \right] dx \nonumber \\
                                           & +  \int_{\Omega} \eta\nabla \chi_2\cdot  \nabla \sigma dx, \\
\label{wf3} B_3^\chi(\chi_3,\phi,\sigma,p) & = \int_{\Omega} \chi_3 \left[ p_t + \gamma_p p -\alpha_h - \alpha_{ch}\phi \right] dx + \int_{\Omega} D\nabla \chi_3\cdot  \nabla p dx.
\end{align}
We also introduce a weak form of the adjoint problem defined by Eqns.~\eqref{34}--\eqref{36}, which is stated as: find $w\in V_0$, $z\in V$, and $q\in V$ such that 
\begin{align} 
\label{b2wf1} & B_1^\psi(\psi_1,w,z,q)=0\quad\text{for all}\quad \psi_1\in V_0, \\
\label{b2wf2} & B_2^\psi(\psi_2,w,z,q)=0\quad\text{for all}\quad \psi_2\in V,   \\
\label{b2wf3} & B_3^\psi(\psi_3,w,z,q)=0\quad\text{for all}\quad \psi_3\in V, 
\end{align}
where,
\begin{align} 
\label{dwf1} 
B_1^\psi(\psi_1,w,z,q)=& - \int_{\Omega} \psi_1 w_t dx + \int_{\Omega} \lambda\nabla \psi_1\cdot  \nabla w dx + \int_{\Omega} \psi_1 \left[ \gamma_{ch} \sigma^\ast + (S^\ast-S_{ch}) \right] z  dx \nonumber\\
& + \int_{\Omega} \psi_1 \left[ F^{\prime\prime}(\phi^\ast) + h^{\prime\prime}(\phi^\ast)(U^\ast-m(\sigma^\ast)) \right] w  dx  -\int_{\Omega} \psi_1\alpha_{ch}q dx \nonumber  \\
& - \int_{\Omega} \psi_1 k_1(\phi^\ast-\phi_Q) dx, \\
\label{dwf2} 
B_2^\psi(\psi_2,w,z,q)=& - \int_{\Omega} \psi_2 z_t dx + \int_{\Omega} \eta\nabla \psi_2\cdot  \nabla z dx + \int_{\Omega} \psi_2 z \left( \gamma_h + \gamma_{ch}\phi^\ast \right)  dx  \nonumber\\
& - \int_{\Omega} \psi_2 w m^\prime(\sigma^\ast)h^\prime(\phi^\ast) dx,  \\
\label{dwf3} 
B_2^\psi(\psi_2,w,z,q)=& -\int_{\Omega} \psi_3 q_t dx +  \int_{\Omega} D\nabla \psi_3\cdot  \nabla q dx + \int_{\Omega} \psi_3 q\gamma_p  dx \nonumber \\
& - \int_{\Omega} \psi_3k_4\left( P^\ast - P_\Omega\right)^{+} dx.
\end{align}
The above-defined weak forms are discretized by defining finite-dimensional spaces $V^h\subset V$ and $V^h_0\subset V_0$. We construct the discrete spaces using  $C^1$-continuous quadratic B-splines \cite{hughes2005isogeometric}. The space $V^h$ is defined as $V^h=\hbox{\rm span}\{N_A\}_{A=1,\dots,n_b}$, where $n_b=\dim(V^h)$ and the $N_A$'s are multivariate splines. We use the superscript $h$ to denote finite-dimensional approximations to the exact solution of the forward and adjoint problems. For example, $\phi^h(t,x)=\sum_{A=1}^{n_b}\phi_A(t)N_A(x)$, where the $\phi_A$'s, called control variables, are the unknowns of the problem. The functions $\sigma^h$, $p^h$, $w^h$, $z^h$, and $q^h$ are defined analogously. The functions that belong to $V^h_0$ will have some control variables constrained to ensure that Dirichlet boundary conditions are satisfied. The Neumann boundary conditions relevant to this problem are naturally enforced  within the weak form.

\subsection{Time integration}

Our time integration scheme is based on the generalized-$\alpha$ algorithm \cite{chung1993time,jansen2000generalized}. Let us call $\boldsymbol{\Phi}$ the global vector of degrees of freedom associated with the unknown $\phi^h$, i.e., $\boldsymbol{\Phi}=\{\phi_A\}_{A=1,\dots,n_b}$. Likewise, we will also make use of the vectors $\boldsymbol{\Sigma}=\{\sigma_A\}_{A=1,\dots,n_b}$ and $\boldsymbol{P}=\{p_A\}_{A=1,\dots,n_b}$. We introduce the residual vector for the forward problem
\begin{equation}
\vec{\rm Res^F} = \{ \boldsymbol{R}^\phi,\boldsymbol{R}^\sigma, \boldsymbol{R}^p\} 
\end{equation}
where $\boldsymbol{R}^\phi=\{R^\phi_A\}_{A=1,\dots,n_b}$, $\boldsymbol{R}^\sigma=\{R^\sigma_A\}_{A=1,\dots,n_b}$, and $\boldsymbol{R}^p=\{R^p_A\}_{A=1,\dots,n_b}$, such that
\begin{align} 
R^\phi_A   & = B_1^\chi(N_A,\boldsymbol{\Phi},\boldsymbol{\Sigma},\boldsymbol{P}),\\
R^\sigma_A & = B_2^\chi(N_A,\boldsymbol{\Phi},\boldsymbol{\Sigma},\boldsymbol{P}),\\
R^p_A      & = B_3^\chi(N_A,\boldsymbol{\Phi},\boldsymbol{\Sigma},\boldsymbol{P}).
\end{align}
We call $\boldsymbol{U}_n=\{\boldsymbol{\Phi}_n,\boldsymbol{\Sigma}_n,\boldsymbol{P}_n\}$ the time-discrete approximation to the control variables of the forward problem at time $t_n$. In the forward problem, we calculate $\boldsymbol{U}_{n+1}$ from $\boldsymbol{U}_n$ by enforcing the equation
\begin{equation}\label{resi}
\vec{\rm Res^F}(\dot{\boldsymbol{U}}_{n+\alpha_m},\boldsymbol{U}_{n+\alpha_f})=\vec 0,
\end{equation}
where
\begin{alignat}{2}
&\vec U_{n+1}              &= &\,\vec U_n+(t_{n+1}-t_n)\dot{\vec U}_n+\gamma(t_{n+1}-t_n)(\dot{\vec U}_{n+1}-\dot{\vec U}_n),  \\
&\dot{\vec U}_{n+\alpha_m} &= &\,\dot{\vec U}_n+\alpha_m(\dot{\vec U}_{n+1}-\dot{\vec U}_n),  \\
&{\vec U}_{n+\alpha_f}     &= &\,{\vec U}_n+\alpha_f({\vec U}_{n+1}-{\vec U}_n). \label{resf}
\end{alignat}
Here, $t_{n+1}-t_n=\Delta t_n>0$ is the time step, and $\alpha_m$, $\alpha_f$, and $\gamma$ are real-valued parameters that define the accuracy and stability of the algorithm.

For the adjoint problem, we define the global vectors of degrees of freedom $\vec W$, $\vec Z$, and $\vec Q$, corresponding, respectively, to the discrete functions $w^h$, $z^h$, and $q^h$. The residual vector for the adjoint problem is 
\begin{equation}\label{resaii}
\vec{\rm Res^A}=\{ \boldsymbol{R}^w,\boldsymbol{R}^z, \boldsymbol{R}^q\} 
\end{equation}
where $\boldsymbol{R}^w=\{R^w_A\}_{A=1,\dots,n_b}$, $\boldsymbol{R}^z=\{R^z_A\}_{A=1,\dots,n_b}$, and $\boldsymbol{R}q=\{R^q_A\}_{A=1,\dots,n_b}$, such that
\begin{align} 
R^w_A  & = B_1^\psi(N_A,\vec W,\vec Z, \vec Q), \\
R^z_A  & = B_2^\psi(N_A,\vec W,\vec Z, \vec Q), \\
R^q_A  & = B_3^\psi(N_A,\vec W,\vec Z, \vec Q).
\end{align}
Let us call $\vec Y_n=\{\vec W_n,\vec Z_n,\vec Q_n\}$ the global vector of degrees of freedom for the adjoint problem at $t=t_n$. The adjoint problem is solved backwards in time starting from data at $t=T$, so for this problem $\Delta t_n = t_n - t_{n+1} < 0$. The equation that allows us to determine $\vec Y_n$ from $\vec Y_{n+1}$ is 
\begin{equation}\label{resai}
\vec{\rm Res^A}(\dot{\boldsymbol{Y}}_{n+\alpha_m},\boldsymbol{Y}_{n+\alpha_f})=\vec 0
\end{equation}
where
\begin{alignat}{2}
&\vec Y_{n}              &= &\,\vec Y_{n+1}+(t_{n}-t_{n+1})\dot{\vec Y}_{n+1}+\gamma(t_{n}-t_{n+1})(\dot{\vec Y}_{n}-\dot{\vec Y}_{n+1}),  \\
&\dot{\vec Y}_{n+\alpha_m} &= &\,\dot{\vec Y}_{n+1}+\alpha_m(\dot{\vec Y}_{n}-\dot{\vec Y}_{n+1}),  \\
&{\vec Y}_{n+\alpha_f}     &= &\,{\vec Y}_{n+1}+\alpha_f({\vec Y}_{n}-{\vec Y}_{n+1}). \label{fres}
\end{alignat}
As shown in \cite{jansen2000generalized}, the generalized-$\alpha$ algorithm can be made $A$-stable and second-order accurate by taking $\rho_\infty\in[0,1]$ and 
\begin{equation}\label{alphas}
\alpha_m=\frac{1}{2}\left(\frac{3-\rho_\infty}{1+\rho_\infty}\right),\quad \alpha_f=\frac{1}{1+\rho_\infty},\quad \gamma=\frac{1}{2}+\alpha_m-\alpha_f.
\end{equation}
All the calculations presented in this paper were performed taking $\rho_\infty=1/2$ and using Eq.~\eqref{alphas}. We linearized the nonlinear algebraic equations defined by \eqref{resi}--\eqref{resf} and \eqref{resai}--\eqref{fres} by using the Newton-Raphson algorithm. The convergence criterion to advance from one time step to the next one was that the individual residuals (e.g., $\vec R^\phi$, $\vec R^\sigma$, and $\vec R^p$ for the forward problem) are reduced to $\epsilon_{NL}$ of its initial value. The linear systems that result after linearization are solved using GMRES \cite{Saad1986} with diagonal preconditioner up to a predefined tolerance $\epsilon_L$ or a maximum number of iterations. 

\subsection{Optimal control algorithm}\label{sda}

We solve numerically the optimal control problem using the steepest-descent gradient method; see, e.g., Ref. \cite{Arnautu} (Algorithm 2.2). As noted in Remark~4.7, we will consider the controls $U$ and $S$ to be exclusively time-dependent in our simulations, and so we build the steepest-descent gradient algorithm accordingly here. We determine the optimal functions $U$ and $S$ constructing a sequence of approximations $\{U_k\}_{k\geq1}$ and $\{S_k\}_{k\geq1}$ with $U_0$ and $S_0$ given. In what follows, we describe our algorithm to calculate ($U_{k+1}, S_{k+1}$) from $(U_{k},\, S_{k})$, which involves 7 main steps:
\begin{enumerate}[leftmargin=.5cm]
\item[] \textbf{Step 1}: Compute the functions $(\phi_{k},\sigma_{k},p_{k})$ solving the forward problem (\ref{1})--(\ref{5}) with $U=U_k$ and $S=S_k$.
\item[] \textbf{Step 2}: Compute the functions $(w_{k},z_{k},q_{k})$ solving the adjoint problem \eqref{34}--\eqref{38} with $(U^\ast,S^\ast)=(U_k,S_k)$ and $(\phi^{\ast},\sigma^{\ast},p^{\ast})=(\phi_{k},\sigma_{k},p_{k})$.
\item[] \textbf{Step 3}: Evaluate the gradient of $J$ at $(U_k,S_k)$ using Eqns.~\eqref{54-1}--\eqref{54-2}, i.e., 
\begin{alignat}{2}
	d_{U_{k}} &=\nabla _{U}J(U_{k},S_{k})=k_{6}U_{k}-\int_{\Omega }h^{\prime}(\phi_{k})w_{k}dx \\
	d_{S_{k}} &=\nabla _{S}J(U_{k},S_{k})=k_{7}S_{k}-\int_{\Omega }\phi_{k}z_{k}dx.
\end{alignat}
\item[] \textbf{Step 4}: Check if any of the following convergence criteria is satisfied:
\begin{enumerate}[leftmargin=.3cm]
	\item[] Criterion 1: 
	\begin{align}
	||d_{U_k}||^{2}_{L^2(0,T)}&<\varepsilon_{SD1}||d_{U_0}||^{2}_{L^2(0,T)}, \notag\\  ||d_{S_k}||^{2}_{L^2(0,T)}&<\varepsilon_{SD1}||d_{S_0}||^{2}_{L^2(0,T)}
	\end{align}	
	\item[] Criterion 2: 
     \begin{align}
     ||d_{U_k}-d_{U_{k-1}}||^{2}_{L^2(0,T)}&<\varepsilon_{SD2}||d_{U_{k-1}}||^{2}_{L^2(0,T)}, \notag\\ ||d_{S_k}-d_{S_{k-1}}||^{2}_{L^2(0,T)}&<\varepsilon_{SD2}||d_{S_{k-1}}||^{2}_{L^2(0,T)}
     \end{align} 
\end{enumerate}
where $\varepsilon_{SD1}$ and $\varepsilon_{SD2}$ are predefined tolerances. If any of the convergence criteria is satisfied, the iterative process to compute the optimal functions ends. If none of the convergence criteria are satisfied, we proceed to Step 5.
\item[] \textbf{Step 5}: Compute a pool of potential values for $(U_{k+1},S_{k+1})$ using the updates
\begin{alignat}{2}
U_{k, \mu_j} &= U_{k} - \mu_j d_{U_{k}} \quad\mbox{for}\quad j=1,\dots,N_\mu \\
S_{k, \mu_j} &= S_{k} - \mu_j d_{S_{k}} \quad\mbox{for}\quad j=1,\dots,N_\mu
\end{alignat}
where $\mu_j=j/N_\mu$ with $j=1,\dots,N_\mu$. Here, $N_\mu>1$ is an integer parameter of the algorithm.
\item[] \textbf{Step 6}: Select the best functions in the pool $\{U_{k+1,\mu_j}\}_{j=1,\dots,N_\mu}$, $\{S_{k+1,\mu_j}\}_{j=1,\dots,N_\mu}$ finding $j^\ast$ such that
\begin{equation}
J(U_{k+1, \mu_{j^\ast}}, S_{k+1, \mu_{j^\ast}})=\min_{j=1,\dots,N_\mu}\{J(U_{k+1, \mu_j}, S_{k+1,\mu_j})\}.
\end{equation}
\item[] \textbf{Step 7}: Define the $(k+1)^{th}$ iteration as
\begin{alignat}{2}
U_{k+1} &= U_{k+1, \mu_{j^\ast}} \\
S_{k+1} &= S_{k+1, \mu_{j^\ast}}.
\end{alignat}
\end{enumerate}
This completes one iteration of the steepest descent algorithm. The process is then repeated until one of the convergence criteria defined in Step 4 is satisfied or a predefined maximum number of iterations is reached.

\section{Simulation study of the optimal control problem}\label{simstudy}

\subsection{Description}\label{Rintro}

In this section, we explore the behavior of the optimal control problem stated in Section~\ref{optcontrol} for time-dependent controls $U(t)$ and $S(t)$. To this end, we carry out a simulation study featuring a prostatic tumor treated with combined cytotoxic and antiangiogenic therapy. Thus, the aim of these simulations is to compute the optimal drug-na\"{i}ve cytotoxic and antiangiogenic effects, $U(t)$ and $S(t)$, respectively, to effectively treat the tumor according to our model of PCa growth. We run this study over the time of a single cycle of combined therapy with docetaxel and bevacizumab, i.e., $T=21$ days \cite{Kelly2012,Mottet2018,CGLMRR-1}. In Section~\ref{drugprotocols} we will use the results of this study to explore alternative therapeutic strategies for specific drugs.

In the following, we present values of the model parameters used in the simulations, the choice of terms and corresponding weighting constants in the objective functional $J$, the construction of initial guesses and the admissible space for $U(t)$ and $S(t)$, and the computational details of the specific implementation of the algorithms presented in Section~\ref{numerics} for this simulation study.

\subsubsection{Parameters of the PCa growth model}
We consider an aggressive case of PCa. This scenario corresponds to a tumor having a high Gleason score, which is a routine clinical variable associated to cancer aggressiveness \cite{Mottet2018}. To model this instance of PCa growth, we select $\rho$ and $A$ in Eq.~(\ref{msigma}) based on the average values of tumor cell proliferation and apoptosis, tumor doubling times, and serum PSA doubling times previously reported in the literature for tumors with high Gleason score \cite{Berges1995,Schmid1993,CGLMRR-1}.

The values of the other parameters in our PCa growth model have been adopted from previous studies \cite{lorenzo2016tissue,lorenzo2017hierarchically,lorenzo2019computer,Xu2016,CGLMRR-1}. Table \ref{parameters} summarizes the parameters participating in Eqns.~\eqref{1}--\eqref{3} and their values in the simulation study presented herein.

\begin{table}[t]
\setlength{\tabcolsep}{0pt}
\caption{List of parameters in Eqns.~\eqref{1}--\eqref{3} and their corresponding values in the simulations presented herein.  }
{\small 
{\begin{tabular}{p{0.4\linewidth}p{0.125\linewidth}p{0.35\linewidth}p{0.125\linewidth}}

\toprule
 Parameter & Notation & Value & Reference\\ 
\midrule

\multicolumn{4}{@{}l}{\textbf{Tumor dynamics}} \\

Diffusivity of the tumor phase field  & 
$\lambda$  & 
640 \si{\micro\metre}$^2$/day & 
\cite{lorenzo2016tissue,Xu2016}\\ 

Tumor mobility &
$M$ & 
2.5 1/day &
\cite{lorenzo2016tissue,Xu2016}\\

Net proliferation scaling factor & 
$m_{ref}$ & 
7.55$\cdot 10^{-2}$ 1/day & 
\cite{Xu2016}\\

Scaling reference for proliferation rate & 
$\bar{K}_\rho$  & 
1.50$\cdot 10^{-2}$ 1/day & 
\cite{Berges1995,Schmid1993}\\

{Proliferation rate}  & 
{$K_\rho$} & 
1.50$\cdot 10^{-2}$  1/day &
\cite{Berges1995,Schmid1993}\\

Scaling reference for apoptosis rate & 
$\bar{K}_A$  & 
2.10$\cdot 10^{-2}$ 1/day & 
\cite{Berges1995,Schmid1993}\\ 
                                              
{Apoptosis rate} & 
{$K_A$} & 
1.37$\cdot 10^{-2}$  1/day &
\cite{Berges1995,Schmid1993}\\

\multicolumn{4}{@{}l}{ } \\ 
\multicolumn{4}{@{}l}{\textbf{Nutrient dynamics }} \\
Nutrient diffusivity & 
$\eta$ & 
$6.4 \cdot 10^{4}$ \si{\micro\metre}$^2$/day & 
\cite{lorenzo2016tissue}\\

Nutrient supply in healthy tissue & 
$S_h$ & 
$2$ g/L/day & 
\cite{lorenzo2016tissue}\\

Nutrient supply in tumor tissue  & 
$S_c$ 
& 2.75  g/L/day &
\cite{lorenzo2016tissue}\\

Nutrient uptake by healthy tissue & 
$\gamma_h$ & 
$2$ g/L/day & 
\cite{CGLMRR-1}\\

Nutrient uptake by tumor tissue  & 
$\gamma_c$ &
17    g/L/day & 
\cite{CGLMRR-1}\\

\multicolumn{4}{@{}l}{ } \\ 
\multicolumn{4}{@{}l}{\textbf{Tissue PSA dynamics}} \\
Tissue PSA diffusivity & 
$D$ &  
640 \si{\micro\metre}$^2$/day & 
\cite{lorenzo2016tissue}\\

Healthy tissue PSA production rate & 
$\alpha_h$ & 
1.712$\cdot 10^{-2}$ ng/mL/cc/day & 
\cite{lorenzo2016tissue}\\

Tumoral tissue PSA production rate & 
$\alpha_c$ & 
$\alpha_c=15\alpha_h$ & 
\cite{lorenzo2016tissue}\\

Tissue PSA natural decay rate & 
$\gamma_p$ & 
0.274 1/day & 
\cite{lorenzo2016tissue}\\

\bottomrule
\end{tabular}}}
\label{parameters}
\end{table}

\subsubsection{Objective functional}\label{objsim}
As explained in Section~\ref{optcontrol}, we selected only some $k_i$, $i=1,\ldots,7$ to be nonzero in our simulations. In particular, we consider three versions of the functional $J$ defined in Eq.~\eqref{J} corresponding to three different sets of nonzero constants $k_i$, $i=1,\ldots,7$, as follows:

\begin{eqnarray}
J_1(U,S) &=& \frac{k_{1}}{2}\int_{Q_{T}}(\phi (t,x)-\phi _{Q})^{2}dxdt 
          + \frac{k_{2}}{2}\int_{\Omega }(\phi (T,x)-\phi _{\Omega })^{2}dx \notag \\  \notag
          &&+ \frac{k_{4}}{2}\int_{0}^{T}\left[ \left(\int_{\Omega }p(t,x)dx-p_{\Omega }(t)\right) ^{+}\right] ^{2}dt \\ 
          &&+\frac{k_{6}}{2}\int_0^T U^{2}(t)dt 
            +\frac{k_{7}}{2}\int_0^T S^{2}(t)dt,  \label{J1}
\end{eqnarray}
\begin{eqnarray}
J_2(U,S) &=& \frac{k_{2}}{2}\int_{\Omega }(\phi (T,x)-\phi _{\Omega })^{2}dx 
          + \frac{k_{4}}{2}\int_{0}^{T}\left[ \left(\int_{\Omega }p(t,x)dx-p_{\Omega }(t)\right) ^{+}\right] ^{2}dt \notag \\
          &&+\frac{k_{6}}{2}\int_0^T U^{2}(t)dt 
            +\frac{k_{7}}{2}\int_0^T S^{2}(t)dt,  \label{J2}
\end{eqnarray}
\noindent and
\begin{eqnarray}
J_3(U,S) &=& k_{3}\int_{\Omega }\phi (T,x)dx 
          + \frac{k_{4}}{2}\int_{0}^{T}\left[ \left(\int_{\Omega }p(t,x)dx-p_{\Omega }(t)\right) ^{+}\right] ^{2}dt \notag \\
          &&+\frac{k_{6}}{2}\int_0^T U^{2}(t)dt 
            +\frac{k_{7}}{2}\int_0^T S^{2}(t)dt.  \label{J3}
\end{eqnarray}
\noindent We choose $\phi_Q=0$ and $\phi_\Omega=0$, such that the optimization problem aims at killing the tumor. Consequently, we set $p_\Omega=\alpha_h \mbox{vol}(\Omega)/\gamma_p$, which is the value of serum PSA that would be obtained if the complete computational domain $\Omega$ were comprised of healthy tissue. We did not consider the term multiplied by $k_5$ in the original definition of $J$ in Eq.~\eqref{J} because, in general, serum PSA can be obtained with reasonable time resolution in experimental and clinical settings \cite{Mottet2018}. 

In this simulation study, we will explore the performance of our optimal control formulation under an array of values of the weighting constants in $J_1$, $J_2$, and $J_3$. In simulations with $J_1$, we impose $k_1=k_2$ due to the similarity of the terms with $k_1$ and $k_2$. Indeed, notice that the term with $k_2$ introduces in the transversality condition for the adjoint variable $w$ (Eq.~(\ref{38})) a similar term to the forcing induced by the $k_1$ term in the equation for $w$ (Eq.~(\ref{34})).
We further choose $k_4=k_1$ such that we use the same weighting for the terms accounting for tumor volume and serum PSA in the objective functional. Likewise, we impose $k_4=k_2$ and $k_4=k_3$ for $J_2$ and $J_3$, respectively. We also fix $k_6=1$ and $k_7=1$ in all simulations in Section~\ref{resJ}, such that we do not favor one therapy over the other within our modeling framework. These choices make $J_1$, $J_2$, and $J_3$ solely dependent on $k_1$, $k_2$, and $k_3$, respectively; which facilitates the analysis of this first simulation study of the optimal control problem presented herein. The specific values of the weighting constants employed in each simulation are reported along with the corresponding results in Section~\ref{resJ}. 

\subsubsection{Initial guess and maximum values of the controls}
Pharmacodynamic studies of common cytotoxic and antiangiogenic drugs for cancer treatment usually show an exponential decay in drug concentration following the systemic delivery of the prescribed dose \cite{Baker2004,Tije2005,Gordon2001,Lu2008}.
Additionally, previous efforts to model cytotoxic and antiangiogenic drug effects usually rely on a linear dependence on the drug concentration \cite{Kohandel2007,Powathil2007,Hinow2009,Jarrett2018,Benzekry2013,Bogdanska2017,CGMR-AMO,GLR,CGLMRR-1}. 
We adopted this approach for the formulation of the initial guess of cytotoxic drug effects on tumor dynamics $U_0(t)$, which is hence given by
\begin{equation}\label{eq_u}
U_0(t)=m_{ref}\beta_c d_c e^{-\frac{t}{\tau_c}},
\end{equation}
\noindent where $\beta_c$ measures the effect of the cytotoxic drug on tumor dynamics per unit of drug dose delivered, $d_c$ is the prescribed dose of the cytotoxic drug, and $\tau_c$ stands for the mean lifetime of the chemotherapeutic drug. Likewise, we define the initial guess of the antiangiogenic effect on nutrient supply $S_0(t)$ as
\begin{equation}\label{eq_v}
S_0(t)=\beta_a d_a e^{-\frac{t}{\tau_a}},
\end{equation}
\noindent where $\beta_a$ measures the effect of the antiangiogenic drug on the nutrient supply per unit of drug dose delivered, $d_a$ is the prescribed dose of antiangiogenic drug, and $\tau_a$ denotes the mean lifetime of the antiangiogenic drug. Notice that this choice of $U_0(t)$ and $S_0(t)$ matches our previous modeling approach \cite{CGLMRR-1}. Additionally, observe that under the assumption of linear dependence of $U(t)$ and $S(t)$ on the cytotoxic and antiangiogenic drug concentrations, respectively, the last two terms in our optimal control problem functional $J$ in Eq.~\eqref{J} are penalizing large concentrations of cytotoxic and  antiangiogenic drugs \cite{Kohandel2007,Powathil2007,Hinow2009,Jarrett2018,Benzekry2013,Bogdanska2017,CGMR-AMO,GLR,CGLMRR-1}.

Docetaxel is considered the gold standard drug in cytotoxic chemotherapy of advanced PCa \cite{Mottet2018,Kelly2012,Eisenberger2012}. Additionally, bevacizumab has been extensively investigated in single-drug antiangiogenic or combined cytotoxic-antiangiogenic drug protocols for advanced PCa \cite{Kelly2012,Antonarakis2012,picus2011phase}. Therefore, we use the standard dosage and pharmacodynamic properties of these two drugs to calculate $U_0(t)$ and $S_0(t)$ in our simulations, respectively. The values of the parameters involved in the dynamics of drug effects in Eqns.~\eqref{eq_u}-\eqref{eq_v} are provided in Table~\ref{u0s0}.

\begin{table}[t]
\setlength{\tabcolsep}{0pt}
\caption{List of parameters in the definition of the initial guess for the controls $U_0(t)$ and $S_0(t)$ given in Eqns.~\eqref{eq_u}-\eqref{eq_v} and their corresponding values. }
{\small 
{\begin{tabular}{p{0.4\linewidth}p{0.125\linewidth}p{0.35\linewidth}p{0.125\linewidth}}

\toprule
 Parameter & Notation & Value & Reference\\ 
\midrule

\multicolumn{4}{@{}l}{\textbf{Cytotoxic chemotherapy}} \\
Mean lifetime of cytotoxic drug & 
$\tau_c$  & 
5 days & 
\cite{Baker2004,Tije2005}\\ 

Cytotoxic drug effect & 
$\beta_c$ & 
1.59$\cdot 10^{-2}$ 1/(mg/m$^2$) &
\cite{CGLMRR-1}\\

Cytotoxic drug dose & 
$d_c$ &  
75 mg/m$^2$ &
\cite{Mottet2018,Kelly2012}\\

\multicolumn{4}{@{}l}{ } \\                                               
\multicolumn{4}{@{}l}{\textbf{Antiangiogenic therapy}} \\
Mean lifetime of antiangiogenic drug & 
$\tau_a$  & 
30 days & 
\cite{Gordon2001,Lu2008}\\ 

Antiangiogenic drug effect & 
$\beta_a$  & 
0.04 g/L/day/(mg/kg) & 
\cite{CGLMRR-1}\\

Antiangiogenic drug dose & 
$d_a$  & 
15 mg/kg & 
\cite{Antonarakis2012,Kelly2012}\\

\bottomrule
\end{tabular}}}
\label{u0s0}
\end{table}

Additionally, the admissible values of the controls $U(t)$ and $S(t)$ are respectively bounded by maximum values $U_{\max}$ and $S_{\max}$ in our formulation of the optimal control problem (see Section~\ref{optcontrol}).
We set $U_{\max}=0.012$ 1/day and $S_{\max}=0.80$ g/L/day. These values respectively correspond to a maximum dose of 100 mg/m$^2$ of docetaxel and 20 mg/kg of bevacizumab under the modeling assumptions of Eqns.~\eqref{eq_u}-\eqref{eq_v} \cite{Hainsworth2004,Engels2005,Montero2005,Kazazi2010,Lyseng2006}.

While we have used docetaxel and bevacizumab to define $U_0(t)$, $S_0(t)$, $U_{\max}$, and $S_{\max}$, we want to remark that the optimal controls $U(t)$ and $S(t)$ describe the optimal cytotoxic and antiangiogenic effects according to our PCa growth model, which may also be obtained with different dosages of other drugs showing different pharmacodynamics. In Section~\ref{drugprotocols}, we provide and example of how to estimate some drug protocols that would approximately yield the optimal $U(t)$ and $S(t)$ obtained in the simulations of our optimal control problem.

\subsection{Computational setup}\label{compsetup}
We  implemented the  numerical algorithms introduced in Section~\ref{numerics} to resolve the forward problem, the adjoint problem, and the steepest descent gradient method by extending our in-house isogeometric codes to simulate PCa growth  \cite{lorenzo2016tissue,lorenzo2017hierarchically,lorenzo2019computer}. These codes were built following the general directions in \cite{hughes2005isogeometric}.

\subsubsection{Space and time discretization}
The computational domain in all simulations is a  square with side length of $L_d=3000$ \si{\micro\metre} and 256 isogeometric elements per side. The time step was set to a constant value of $\Delta t_n=0.1$ days for the forward problem and $\Delta t_n=-0.1$ days for the adjoint problem. 

\subsubsection{Convergence of numerical algorithms}
The convergence of the Newton-Raphson method was set to tolerance $\varepsilon_{NL}=10^{-3}$, while for the GMRES algorithm was set to $\varepsilon_{L}=10^{-3}$ or a maximum of 500 iterations. The convergence for the steepest-descent gradient algorithm was set to $\varepsilon_{SD1}=\varepsilon_{SD2}=10^{-6}$ or a maximum of 100 iterations. To find the update for the optimal $U(t)$ and $S(t)$ in each step of the steepest-descent gradient algorithm, we evaluated $N_\mu=10$ evenly-spaced values of $\mu$ in $(0,1]$, i.e., $\mu_j=j/10$ for $j=1,\ldots,10$ (see Section~\ref{sda}).

\subsubsection{Initial conditions of the forward problem}
We approximate the initial tumor phase field as an ellipsoidal tumor placed in the center of the domain with semiaxes $a=150$ \si{\micro\metre} and $b=200$ \si{\micro\metre} parallel to the domain sides. We implement this initial condition by $L^2$-projecting the hyperbolic tangent function
\begin{equation}
\phi_0(x)=\phi_0(x_1,x_2) = 0.5 - 0.5\tanh\left( 10\left(\sqrt{\frac{(x_1-L_d/2)^2}{a^2} + \frac{(x_2-L_d/2)^2}{b^2}} - 1\right) \right)
\end{equation}
over the quadratic B-spline space supporting our spatial discretization. This operation yields the control variables $\phi_{0,A}=\phi_{A}(0)$, $A=1,\ldots,n_b$, for the spline representation of the phase-field initial condition, i.e., $\phi^h_0(x)=\phi^h(0,x)=\sum_{A=1}^{n_b}\phi_{0,A}N_A(x)$ (see Section \ref{space_disc}).

The initial conditions for the nutrient and the tissue PSA are estimated from $\phi_0$ as
\begin{equation}
\sigma_0 = c^0_\sigma + c^1_\sigma\phi_0
\end{equation}
and
\begin{equation}
p_0 = c^0_p + c^1_p\phi_0.
\end{equation}
The constants $c^0_\sigma$, $c^1_\sigma$, $c^0_p$, and $c^1_p$ are computationally estimated \cite{lorenzo2017hierarchically}, such that $\sigma_0$ and $p_0$ represent a constant value of the nutrient and tissue PSA within the tumor and the host tissue. Hence, we choose $c^0_\sigma=1$ g/L, $c^1_\sigma=-0.8$ g/L, $c^0_p=0.0625$ ng/mL/cc, and $c^1_p=0.7975$ ng/mL/cc.

In all simulations, we initially let the tumor grow untreated for 60 days, i.e., $U(t)=0$ and $S(t)=0$. This enables us to have a reference for the dynamics of the untreated tumor and obtain a good estimate of fields $\phi$, $\sigma$, and $p$ according to our PCa growth model away from the estimated initial conditions \cite{CGLMRR-1}. Then, we proceed to add cytotoxic and antiangiogenic drug effects. To facilitate the ensuing implementation of the steepest-descent algorithm, we reset $t=0$ at the date of drug delivery and we use the values of $\phi$, $\sigma$, and $p$ obtained at the end of the untreated growth phase as initial conditions for the forward problem within our optimal control framework. Hence, the forward and dual problems are run between $t=0$ and $t=T=21$ days,  corresponding to the common duration of a combined therapy cycle for PCa \cite{Kelly2012,Mottet2018,CGLMRR-1}.

\subsection{Simulation results of the optimal control problem}\label{resJ}

\begin{figure}
\centerline{\includegraphics[width=\linewidth]{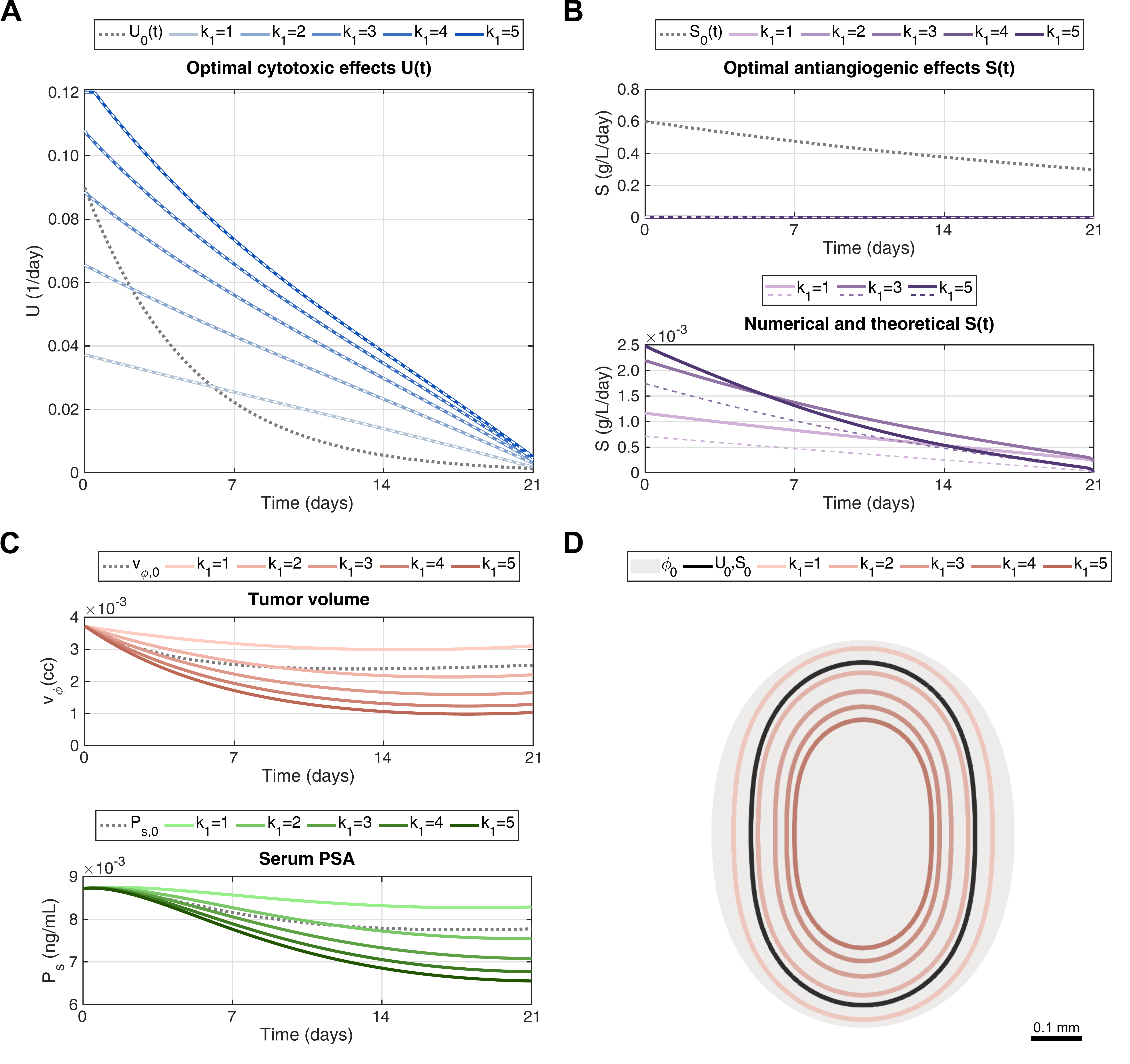}}
\vspace*{8pt}
\caption{Results of the optimal control problem using functional $J_1$ for $k_1=1,2,3,4,5$. (A) Optimal cytotoxic effects $U(t)$ obtained for each value of $k_1$ compared to standard docetaxel therapy ($U_0(t)$, gray dotted line) and the corresponding theoretical estimates calculated with Eq.~(\ref{40-0}) (overlapping white dashed lines). (B) Top: optimal antiangiogenic effects $S(t)$ obtained for each value of $k_1$ compared to standard bevacizumab therapy ($S_0(t)$, gray dotted line) and the corresponding theoretical estimates calculated with Eq.~(\ref{40-0}) (overlapping white dashed lines). Bottom: Detail of optimal $S(t)$ for $k_1=1,3,5$ calculated numerically (solid lines) and with Eq.~(\ref{40-0}) (dashed lines). (C) Evolution of tumor volume $v_\phi$ (top) and serum PSA $P_s$ (bottom) using the optimal $U(t)$ and $S(t)$ obtained via simulation for each value of $k_1$ compared to standard combined therapy, i.e., $v_{\phi,0}$ and $P_{s,0}$  respectively (gray dotted lines). (D) Tumor contours obtained at $t=T$ using the optimal $U(t)$ and $S(t)$ obtained via simulation for each value of $k_1$. The gray area denoted with $\phi_0$ is the tumor region at the onset of the optimal control problem ($t=0$). The black line is the tumor contour at $t=T$ for the standard combined therapy ($U_0$, $S_0$).  }
\label{aggrtJ1}
\end{figure}

\begin{figure}
\centerline{\includegraphics[width=\linewidth]{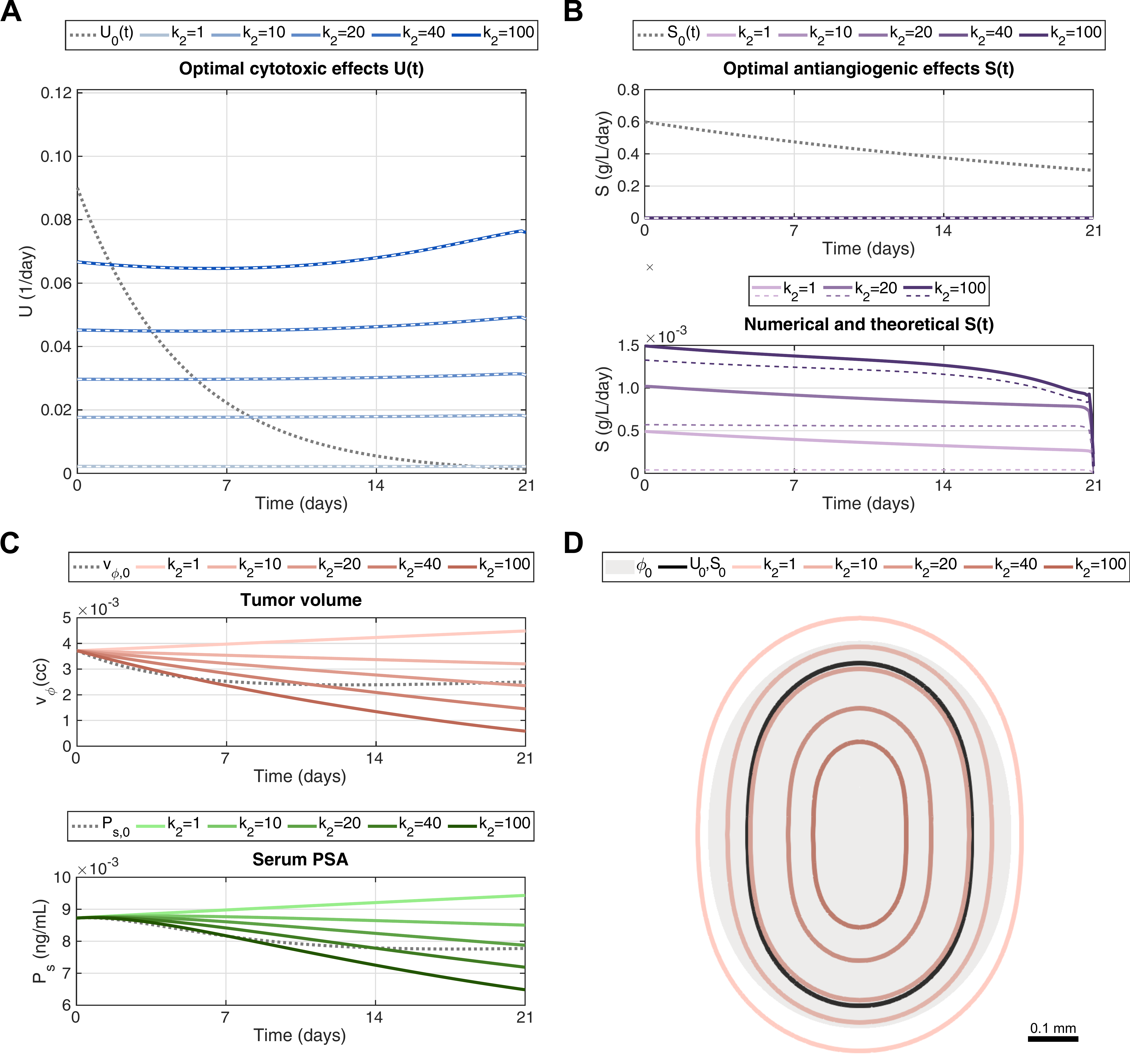}}
\vspace*{8pt}
\caption{Results of the optimal control problem using functional $J_2$ for $k_2=1,10,20,40,100$. (A) Optimal cytotoxic effects $U(t)$ obtained for each value of $k_2$ compared to standard docetaxel therapy ($U_0(t)$, gray dotted line) and the corresponding theoretical estimates calculated with Eq.~(\ref{40-0}) (overlapping white dashed lines). (B) Top: optimal antiangiogenic effects $S(t)$ obtained for each value of $k_2$ compared to standard bevacizumab therapy ($S_0(t)$, gray dotted line) and the corresponding theoretical estimates calculated with Eq.~(\ref{40-0}) (overlapping white dashed lines). Bottom: Detail of optimal $S(t)$ for $k_2=1,20,100$ calculated numerically (solid lines) and with Eq.~(\ref{40-0}) (dashed lines). (C) Evolution of tumor volume $v_\phi$ (top) and serum PSA $P_s$ (bottom) using the optimal $U(t)$ and $S(t)$ obtained via simulation for each value of $k_2$ compared to standard combined therapy, i.e., $v_{\phi,0}$ and $P_{s,0}$ respectively (gray dotted lines). (D) Tumor contours obtained at $t=T$ using the optimal $U(t)$ and $S(t)$ obtained via simulation for each value of $k_2$. The gray area denoted with $\phi_0$ is the tumor region at the onset of the optimal control problem ($t=0$). The black line is the tumor contour at $t=T$ for the standard combined therapy ($U_0$, $S_0$).  }
\label{aggrtJ2}
\end{figure}

\begin{figure}
\centerline{\includegraphics[width=\linewidth]{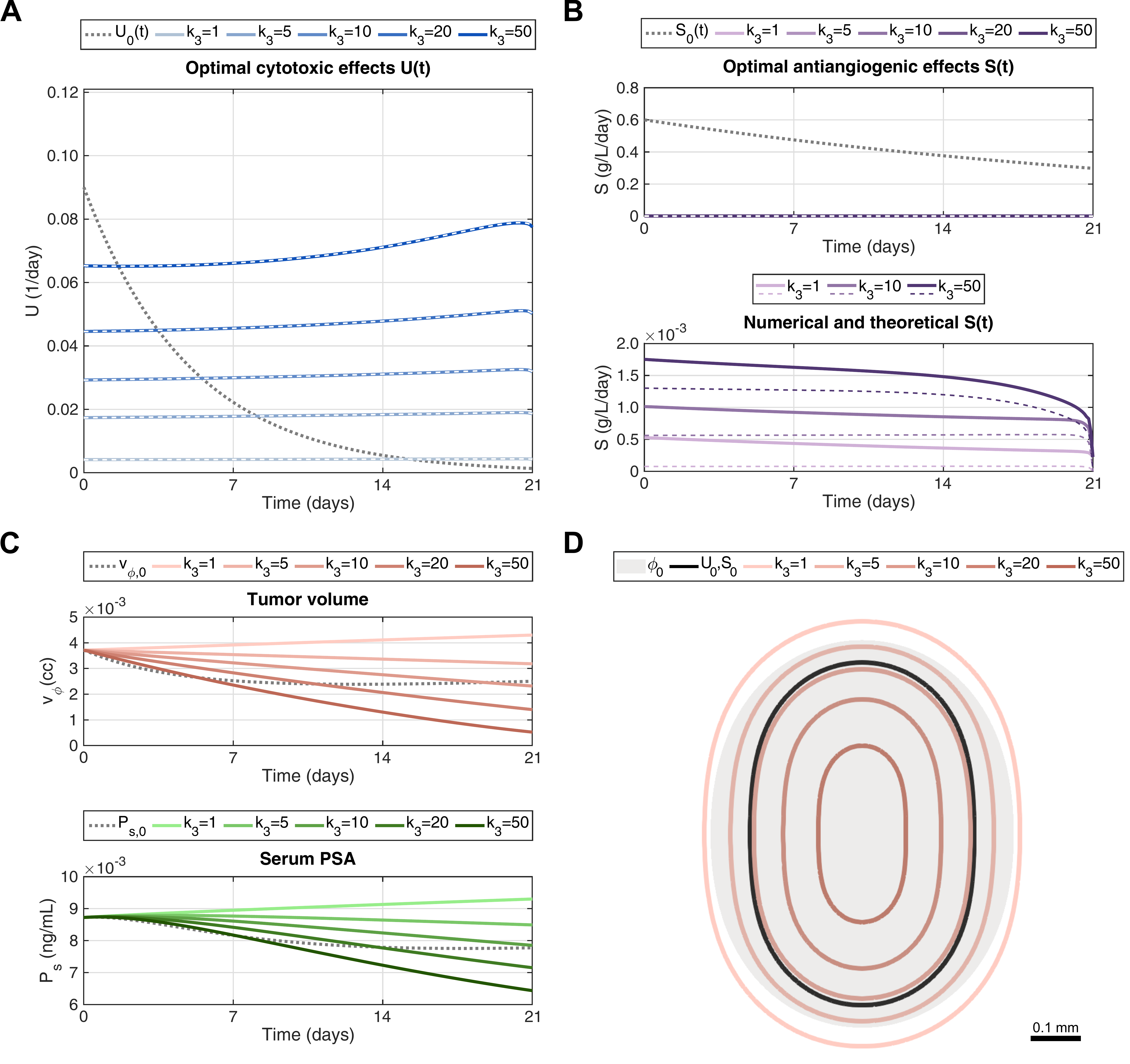}}
\vspace*{8pt}
\caption{Results of the optimal control problem using functional $J_3$ for $k_3=1,5,10,20,50$. (A) Optimal cytotoxic effects $U(t)$ obtained for each value of $k_3$ compared to standard docetaxel therapy ($U_0(t)$, gray dotted line) and the corresponding theoretical estimates calculated with Eq.~(\ref{40-0}) (overlapping white dashed lines). (B) Top: optimal antiangiogenic effects $S(t)$ obtained for each value of $k_3$ compared to standard bevacizumab therapy ($S_0(t)$, gray dotted line) and the corresponding theoretical estimates calculated with Eq.~(\ref{40-0}) (overlapping white dashed lines). Bottom: Detail of optimal $S(t)$ for $k_3=1,10,50$ calculated numerically (solid lines) and with Eq.~(\ref{40-0}) (dashed lines). (C) Evolution of tumor volume $v_\phi$ (top) and serum PSA $P_s$ (bottom) using the optimal $U(t)$ and $S(t)$ obtained via simulation for each value of $k_3$ compared to standard combined therapy, i.e., $v_{\phi,0}$ and $P_{s,0}$ respectively (gray dotted lines). (D) Tumor contours obtained at $t=T$ using the optimal $U(t)$ and $S(t)$ obtained via simulation for each value of $k_3$. The gray area denoted with $\phi_0$ is the tumor region at the onset of the optimal control problem ($t=0$). The black line is the tumor contour at $t=T$ for the standard combined therapy ($U_0$, $S_0$).  }
\label{aggrtJ3}
\end{figure}

Figures~\ref{aggrtJ1}--\ref{aggrtJ3} respectively show the optimal $U(t)$ and $S(t)$ distributions that we obtained in simulations of our optimal control problem using the objective functionals $J_1$, $J_2$, and $J_3$ for an array of values of $k_1$, $k_2$, and $k_3$. These figures also depict the evolution of tumor volume and serum PSA obtained for each optimal solution, along with the tumor contours corresponding to the isosurface $\phi=0.5$ at $t=T$. Larger values of the driving weighting constant $k_i$, $i=1,2,3$, in each corresponding objective functional $J_i$, $i=1,2,3$, increasingly penalize tumor volume and serum PSA (see Section~\ref{objsim}). Consequently, the minimization of $J_i$, $i=1,2,3$, will require at least one of the controls to also increase for larger choices of $k_i$, $i=1,2,3$, because this provides a more intense inhibitory effect on tumor dynamics that ultimately leads to smaller tumor volume and lower serum PSA, as shown in Figs.~\ref{aggrtJ1}--\ref{aggrtJ3}.

In all cases run with $J_1$, the optimal therapy obtained in our simulations involves only a monotonically decreasing cytotoxic effect and virtually no antiangiogenic effect, i.e., $S(t)\approx 0$. The optimal cytotoxic effect even saturates to $U_{\max}$ at early times for $k_1=5$, whereas $S(t)$ still takes negligible values. Thus, these results suggest that strong cytotoxic effects suffice to optimally control the volume and serum PSA of the simulated tumor according to our PCa model. A similar result is obtained for the simulations involving $J_2$ and $J_3$ depicted in Figs.~\ref{aggrtJ2}--\ref{aggrtJ3}. However, the optimal cytotoxic effects $U(t)$ obtained for functionals $J_2$ and $J_3$ are approximately constant for early times and increase in time towards $t=T$, especially for large values of $k_2$ and $k_3$. As the tumor phase field $\phi$ drives serum PSA dynamics in our PCa model (see Section~\ref{pcamodel}), the tumor volume becomes the most powerful clinical quantity of interest participating in the objective functional $J$. Therefore, the increase in optimal $U(t)$ by $t=T$ for larger $k_2$ and $k_3$ is a consequence of $J_2$ and $J_3$ controlling for tumor volume only at $t=T$. 

The simulations presented in Figs.~\ref{aggrtJ1}--\ref{aggrtJ3} show that our optimal control problem enables the calculation of an optimal $U(t)$ providing a similar control of tumor volume and serum PSA as the standard combined therapy with docetaxel and bevacizumab, i.e., $U_0(t)$ and $S_0(t)$ respectively. Figures~\ref{aggrtJ1}--\ref{aggrtJ3} also show that we can obtain optimal $U(t)$ and $S(t)$ exhibiting an improved therapeutic performance over the standard combined protocol as we increase the driving weighting constant $k_i$, $i=1,2,3$. However, we notice that the optimal $U(t)$ obtained with $J_2$ and $J_3$ would require a prolonged exposition to the cytotoxic drug. This would entail a high toxicity that may eventually provoke serious side-effects and may potentially lead to abort the treatment prematurely. Because $J_1$ accounts for tumor volume for all times, we believe that this objective functional enables the computation of superior optimal drug effects precisely tuning treatment to tumor dynamics. We note that the curves of optimal cytotoxic effects $U(t)$ shown in Fig.~\ref{aggrtJ1} do not resemble the exponential trends usually caused by the pharmacodynamics of usual cytotoxic drugs, e.g., the docetaxel curve corresponding to $U_0(t)$ in Figs.~\ref{aggrtJ1}--\ref{aggrtJ3}. In Section~\ref{drugprotocols}, we illustrate how these optimal $U(t)$ time distributions can be leveraged to design novel drug protocols.

Additionally, Figures~\ref{aggrtJ1}--\ref{aggrtJ3} show that the optimal $U(t)$ found in our simulations matches the corresponding theoretical formulation in Eq.~\eqref{40-0} for each value of the driving weighting constant $k_i$, $i=1,2,3$, in $J_i$, $i=1,2,3$, respectively.  
Figures~\ref{aggrtJ1}--\ref{aggrtJ3} also show a detail of optimal $S(t)$ distributions computed numerically and using the theoretical formulation from Eq.~\eqref{40-0} for some simulations. In general, these two optimal $S(t)$ distributions differ in a value that would not be captured by the tolerance of our implementation of the steepest-descent gradient algorithm (see Section~\ref{numerics}). However, the numerical and theoretical solutions for $k_1=5$ virtually coincide. We also observe increasingly better agreement between both distributions for larger values of the driving weighting constant $k_i$, $i=1,2,3$, respectively in $J_i$, $i=1,2,3$. Hence, further reducing the tolerance in the steepest-descent gradient algorithm would make both distributions of optimal $S(t)$ converge. Nevertheless, we want to remark that both the theoretical and numerical optimal coincide for therapeutic purposes here, as both expressions of $S(t)\approx 0$  regardless of the objective functional considered. 

Finally, Appendix B shows ancillary simulations in which we decrease $k_7$ by several orders of magnitude using $J_i$, $i=1,2,3$.
For each of the objective functionals, simulations with $k_7=0.01$ and $k_7=0.001$ render a non-negligible optimal $S(t)$ that matches the corresponding theoretical estimates provided by Eq.~\eqref{40-0}. 
For the simulations with $J_1$ these optimal $S(t)$ solutions show a decreasing profile in time, while for $J_2$ and $J_3$ the corresponding optimal $S(t)$ solutions are mostly constant or slightly increasing towards $t=T$.
The optimal cytotoxic effects $U(t)$ remain practically constant in all simulations for each choice of the objective functional as we vary $k_7$ and the time distribution is virtually the same as the corresponding simulations shown in Figs.~\ref{aggrtJ1}--\ref{aggrtJ3}. 
However, the increase of optimal $S(t)$ has limited effect on tumor volume and serum PSA compared to the parallel increase observed in optimal $U(t)$ in Figs.~\ref{aggrtJ1}--\ref{aggrtJ3}. This suggests that the choices of $k_7$ in Appendix B might artificially bias our optimal control problem and produce unrealistic solutions of limited therapeutic use. 
As we discuss in Section~\ref{disc}, a rigorous analysis of the values of the weighting constants $k_i$, $i=1,\ldots,7$, is required to identify choices representing clinically relevant scenarios for our optimal control framework.

\section{Design of drug protocols inspired by the computed optimal effects}\label{drugprotocols}

\subsection{Methodology}

To illustrate the applicability of the results of our optimal control problem, we calculate an array of model-inspired drug protocols aiming at producing the optimal cytotoxic effect $U(t)$ obtained in one of the simulations of our optimal control problem reported in Section~\ref{simstudy}. In this section, we assume $S(t)=0$, in accordance with the results showed in the previous one. 

We consider four different drug protocols, which we model leveraging the same paradigm that we employed to define $U_0(t)$ in Eq.~\eqref{eq_u}. We consider a single dose of docetaxel, for which we simply calculate $d_c$ in Eq.~\eqref{eq_u}. Then, we consider a new design drug for which we compute $\tau_c$ and the docetaxel equivalent dose $d_c$, such that we can use the same $\beta_c$ in Eq.~\eqref{eq_u}. Finally, we will consider a three-dose protocol of either docetaxel or the new design drug. In this case, we extend the formulation of Eq.~\eqref{eq_u} to
\begin{equation}\label{eq_u3}
U(t)=\sum_{i=1}^3 m_{ref}\beta_c d_{c,i} e^{-\frac{t-t_{c,i}}{\tau_c}}\mathcal{H}(t-t_{c,i}),
\end{equation}
where $\mathcal{H}(t)$ is the Heaviside function. We compute the three doses $d_{c,i}$, $i=1,2,3$ for either drug, together with $\tau_c$ for the new design drug. Additionally, while we fix the delivery of the first dose at $t_{c,1}=0$, we calculate the delivery times of the second and third doses for either drug as well (i.e., $t_{c,2}$ and $t_{c,3}$, respectively).

We compute each of these drug protocols by running a nonlinear least-square fit of the optimal $U(t)$ obtained using $J_1(U,S)$ and $k_1=k_2=k_4=2$, which is shown in Fig.~\ref{aggrtJ1}. We use the trust-region method as provided by the Curve Fitting Toolbox in MATLAB (Release R2020a, The Mathworks, Inc., Natick, Massachusetts, US). The starting dose is set at $d_c=75$ mg/m$^2$ in single-dose protocols. In 3-dose protocols the starting point is three equal doses $d_c=25$ mg/m$^2$ delivered weekly, i.e., with $t_{c,2}$=7 days and $t_{c,3}$=14 days. The starting value for $\tau_c$ for the new design drug is 5 days, corresponding to docetaxel (see Table~\ref{u0s0}). The admissible value ranges are $[0,100]$ mg/m$^2$ for drug doses, $[0,21]$ for the delivery times of the second and third doses in 3-dose protocols, and $[1,20]$ for the characteristic decay time of the new design drug. Additionally, we run a simulation of the forward problem with the fitted drug protocols using the computational setup described in Section~\ref{compsetup} and compare the dynamics of tumor volume and serum PSA with those provided by the optimal $U(t)$. We assess the goodness of fit by means of $R^2$ and the root mean squared error (RMSE).

\subsection{Results}

\begin{figure}[t]
\centerline{\includegraphics[width=\linewidth]{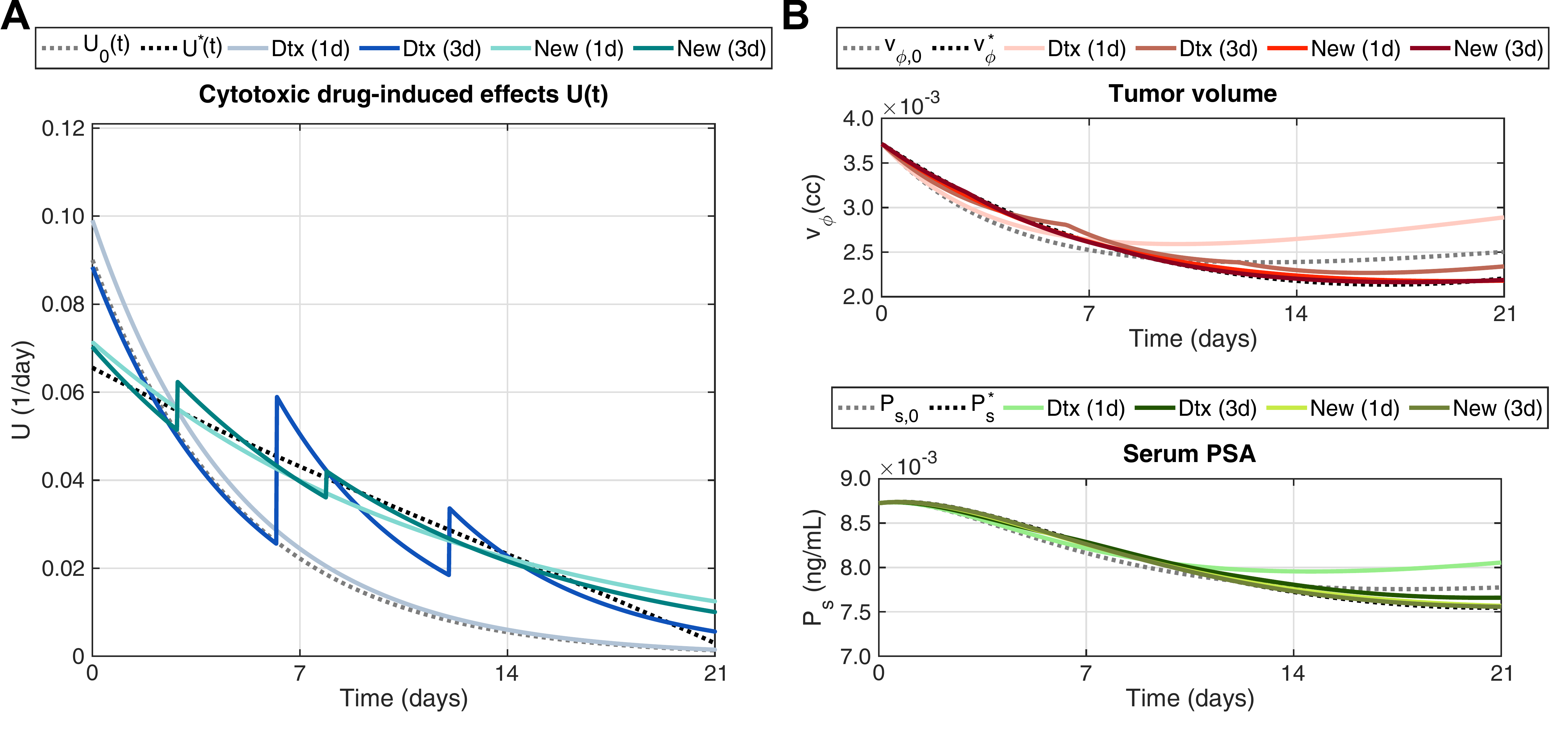}}
\vspace*{8pt}
\caption{Design of drug protocols with docetaxel (Dtx) or a new design drug (New) aiming at reproducing the optimal cytotoxic effect $U^*(t)$ as obtained from the simulation of our optimal control problem using $J_1(U,S)$ and $k_1=k_2=k_4=2$. (A) Cytotoxic drug-induced effects for each of the calculated protocols compared to those of the standard docetaxel protocol $U_0(t)$ and the target optimal cytotoxic effects obtained in our simulations of the optimal control problem $U^*(t)$. (B) Evolution of tumor volume $v_\phi$ and serum PSA $P_s$ for each calculated drug protocol, the standard docetaxel protocol ($v_{\phi,0}$, $P_{s,0}$), and the optimal control problem solution ($v^\ast_{\phi}$, $P^\ast_{s}$). }
\label{drugtests}
\end{figure}

\begin{table}[t]
\setlength{\tabcolsep}{0pt}
\caption{Parameter values obtained in the nonlinear fitting with the trust-region method for each investigated drug protocol.}
{\small 
{\begin{tabular}{p{0.19\linewidth}p{0.135\linewidth}p{0.135\linewidth}p{0.135\linewidth}p{0.135\linewidth}p{0.135\linewidth}p{0.135\linewidth}}

\toprule
Parameters & $d_{c,1}$  & $d_{c,2}$  & $d_{c,3}$  &  $t_{c,2}$  & $t_{c,3}$  & $\tau_c$  \\
(Units) & (mg/m$^2$) & (mg/m$^2$) & (mg/m$^2$) & (day) & (day) & (day) \\

\midrule
1-dose docetaxel & 82.53 & - & - & -  & -  & - \\
3-dose docetaxel & 73.69 & 27.97 & 12.72 & 6.20 & 12.04 & - \\
1-dose new drug  & 59.45 & - & -  & - & - & 12.04\\
3-dose new drug  & 58.49 & 9.20 & 5.03 & 2.85 & 7.90 & 9.16\\
\bottomrule
\end{tabular}}}
\label{fitresults}
\end{table}

\begin{table}[t]
\setlength{\tabcolsep}{0pt}
\caption{Fit quality of optimal cytotoxic effects, tumor volume control, and serum PSA control for each investigated drug protocol.}
{\small 
{\begin{tabular}{p{0.19\linewidth}p{0.135\linewidth}p{0.135\linewidth}p{0.135\linewidth}p{0.135\linewidth}p{0.135\linewidth}p{0.135\linewidth}}

\toprule
           & \multicolumn{2}{@{}l}{\textbf{Cytotoxic effects}} & \multicolumn{2}{@{}l}{\textbf{Tumor volume}} & \multicolumn{2}{@{}l}{\textbf{Serum PSA}}\\
\cmidrule{2-7}
 Statistic & $R^2$ & RMSE    & $R^2$ & RMSE  & $R^2$ & RMSE  \\
 (Units)   & (-)   & (1/day) & (-)   & (cc)  & (-)   & (ng/mL)\\

\midrule
1-dose docetaxel & 0.216 & $1.58\cdot 10^{-2}$ & 0.250 & $3.95\cdot 10^{-4}$ & 0.676 & $2.40\cdot 10^{-4}$ \\
3-dose docetaxel & 0.786 & $8.32\cdot 10^{-3}$ & 0.950 & $1.03\cdot 10^{-4}$ & 0.975 & $6.70\cdot 10^{-5}$ \\
1-dose new drug  & 0.959 & $3.62\cdot 10^{-3}$ & 0.994 & $3.63\cdot 10^{-5}$ & 0.997 & $2.40\cdot 10^{-5}$ \\
3-dose new drug  & 0.976 & $2.78\cdot 10^{-3}$ & 0.999 & $1.69\cdot 10^{-5}$ & 0.999 & $1.02\cdot 10^{-5}$ \\
\bottomrule
\end{tabular}}}
\label{fitgof}
\end{table}

Table~\ref{fitresults} shows the result of the nonlinear fitting of the parameters participating in the formulation of each protocol. Figure~\ref{drugtests} shows the cytotoxic effects generated by each of the investigated protocols, as well as the tumor volume and serum PSA evolution that they would produce according to the corresponding forward simulations of our PCa growth model. Table~\ref{fitgof} shows the fit quality of the cytotoxic effects, tumor volume, and serum PSA for each protocol with respect to those provided by the optimal $U(t)$. 
The single-dose docetaxel protocol uses a slightly larger dose than the standard protocol (75 mg/m$^2$, see Table~\ref{u0s0}), but, in absence of the antiangiogenic drug, the corresponding control of tumor volume and serum PSA is poorer than in the standard protocol. The single-dose docetaxel protocol also produced the worse fitting of the optimal cytotoxic effects. Conversely, the 3-dose docetaxel protocol produces an enhanced fitting of the optimal $U(t)$ and a good control of tumor volume and serum PSA. The best fitting of optimal cytotoxic effects and tumor control was provided by the protocols involving a new design drug. Our calculations suggest that this design drug should have a larger $\tau_c$, approximately twice the corresponding value for docetaxel. This results in slower decay and hence a prolonged cytotoxic effect in time in comparison to docetaxel. We also observe that this new drug would require lower dosage for both the single-dose and 3-dose cases compared to the corresponding protocols using docetaxel. Since the control of tumor volume and serum PSA is virtually the same with the single-dose and 3-dose protocols using the new design drug, we believe that the single-drug protocol would be more suitable for clinical implementation as it would require less visits of the patient to the clinic.

\section{Conclusions}\label{disc}
We present an optimal control theory enabling the calculation of the optimal drug-na\"{i}ve cytotoxic and antiangiogenic effects aiming at successfully treating advanced PCa, which can then be used to design optimal drug compounds and delivery plans. We build this framework departing from our model of PCa growth with cytotoxic and antiangiogenic drug therapies \cite{CGLMRR-1}. We define an objective functional featuring diverse measurements of tumor morphology, tumor volume, and serum PSA, which are common quantities of interest in experimental and clinical studies of advanced PCa.
The functional contains terms accounting for the differences between the computed solution for the prostate tumor or PSA and prescribed values of them, which should be reached during the simulated time or at a final endpoint by an optimal treatment procedure.
We also include penalty terms for cytotoxic and antingiogenic effects, which could be ultimately linked to the corresponding drug concentrations following the usual linear paradigm in the literature \cite{Kohandel2007,Powathil2007,Hinow2009,Jarrett2018,Benzekry2013,Bogdanska2017,CGMR-AMO,GLR,CGLMRR-1}.
From the mathematical viewpoint we prove that the functional reaches its minimum at the optimal values of the cytotoxic and antiangiogenic effects, that is, it provides the optimal design of the treatment which may possibly lead to a final desired result regarding the decrease of the tumor and/or of the PSA. The necessary conditions that should be satisfied by the optimal values are determined by means of the adjoint problem.

In this work, we also propose an algorithm implementing the steepest-descent gradient method to solve our optimal control problem \cite{Arnautu}, which relies on IGA \cite{hughes2005isogeometric} and the generalized-$\alpha$ method \cite{chung1993time,jansen2000generalized} to discretize the involved forward and adjoint problems in space and time, respectively. We then use this algorithm in a simulation study to analyze the performance of our optimal control problem. 
In general, we observed a remarkable agreement between our numerical approximations of the optimal drug effects and the theoretical formulations found in Section~\ref{optcon}. 
Our results show that accounting for the spatio-temporal evolution of the tumor morphology in the objective functional (i.e., $J_1$ in Section~\ref{simstudy}) led to superior optimal drug effects, which precisely adapt the optimal therapy to tumor dynamics and would lead to drug protocols with expected lower toxicity compared to the optimal drug effects obtained by only accounting for the final tumor morphology or volume at $t=T$ (i.e., $J_2$ and $J_3$ in Section~\ref{simstudy} respectively). A previous simulation study on our PCa model \cite{CGLMRR-1} revealed that drug-induced changes in tumor morphology during chemotherapy may contribute to chemoresistance, which supports the inclusion of spatio-temporal dynamics of tumor morphology in the objective functional to refine the design of optimal drug therapies. 

The simulated cases of our optimal control problem in Section~\ref{simstudy} feature $S(t)\approx 0$. We think that these results might align with the studies suggesting that antiangiogenic therapy might not be optimal, and hence that a cytotoxic monotherapy would suffice to optimally treat advanced PCa \cite{Antonarakis2012,Seruga2011,Small2012,mukherji2013angiogenesis,Kelly2012,Tannock2013,Petrylak2015,Smith2016,Sternberg2016}. 
However, several methodological shortcomings have been detected in clinical studies of antiangiogenic therapies, requiring a closer monitoring of pathological events to account for a more direct control of tumor evolution than survival endpoints and serum PSA \cite{Antonarakis2012,Seruga2011,Small2012,mukherji2013angiogenesis}. Some of these monitoring variables are already considered or could be added to the objective functional in our optimal control framework. 
Furthermore, we plan to study the selection of the values of the weighting constants $k_i$, $i=1,\ldots,7$, in the objective functional $J$ to define a constant set enabling the consistent design of personalized optimal drug protocols. We believe that such a selection requires to account for the intrinsic scales in our PCa growth model. Consequently, this update may result in optimal solutions with relevant antiangiogenic effects (e.g.,  Appendix B shows some simulations with non-zero optimal $S(t)$ for $0<k_7\ll 1$).

Additionally, our PCa model only includes antiangiogenic effects as a decrease in nutrient supply. We plan to overcome this simplifying limitation by extending our model to account for the dynamics of the evolving tumor-supporting vascular network \cite{Hanahan2011}. This feature would enable a superior description of nutrient and drug supply to the tumor and provide a clear target for antiangiogenic therapy \cite{Kohandel2007,Powathil2007,Hahnfeldt1999,Hahnfeldt2003,Benzekry2013,Hinow2009,Kremheller2019}.
Accounting for the vascular network would also enable the study of the synergistic or antagonistic effect of combining cytotoxic and antiangiogenic therapies, e.g., whether the reduction of tumor microvasculature caused by antiangiogenic therapy may hamper the effective delivery of cytotoxic drugs to the tumor or whether normalizing the tumor microvasculature would effectively improve the supply of cytotoxic drugs \cite{Jain2005,Mpekris2017}. 
The dynamics of the tumor-induced vascular network can be modeled following different approaches \cite{Vilanova2017}. Hybrid models can capture the evolving morphology of tumor-induced microvasculature with great detail, but usually require intensive computational resources \cite{Frieboes2010,Xu2016,Kremheller2019,Phillips2020,Vavourakis2017,Lima2014}. Conversely, the local density of the vascular network can be modeled using  a continuous  formulation \cite{Kohandel2007,Hahnfeldt1999,Benzekry2013,Vilanova2017,Hormuth2019,Swanson2011}, which would dramatically reduce the computational cost of simulations, especially at tissue and organ scale. This advantage comes at the cost of losing geometrical precision, but it would still enable to account for the tumor-induced vasculature observable in current magnetic resonance imaging  \cite{Hormuth2019,Wu2020}. 

Our model of PCa growth could also be refined by incorporating the effect of tumor-induced mechanical stresses on tumor dynamics, which has been shown to improve tumor forecasting \cite{Weis2015,Lima2016}. Indeed, a recent computational study suggests that the mechanical stresses created by PCa and coexisting benign prostatic hyperplasia may obstruct tumor growth \cite{lorenzo2019computer}. Accounting for tumor-induced mechanical deformation can also improve the modeling of the tumor-supporting vascular network \cite{Hormuth2019,Santos-Oliveira2015}. Furthermore, a poroelastic formulation would enrich the description of nutrient and drug dynamics in the tumor region by taking into account how these phenomena are affected by local changes in mechanical stress and fluid pressure \cite{Fraldi2018,Jain2014,Kremheller2019}. This approach could provide new insights on drug delivery and action in the complex tumor environment, enabling us to refine therapeutic strategies accordingly. Our PCa growth model may also include multiphase formulations accounting for several tumor species with varying responses to the prescribed therapy \cite{Hahnfeldt2003,Jackson2000,Gallaher2018,Wang2016,Lima2014,Rocha2018}. 
To address the computational challenges of these model extensions and hence rationalize the computational resource demand, our numerical methods could be accelerated by implementing adaptive time stepping \cite{Gomez2011,Vilanova2013}, dynamic local adaptivity in the spatial discretization \cite{lorenzo2017hierarchically,Carraturo2019}, and superior algorithms to solve the optimal control problem \cite{Arnautu,Hormuth2018}. 

Finally, we have also presented a strategy that decouples the problem of finding optimal treatment solutions in two phases: first, we compute the optimal therapeutic effect in a drug-independent optimal control framework, and, second, we calibrate alternative drug protocols for specific drug dosage, effects, and pharmacodynamics targeting the calculated optimal therapeutic effect.
The application of this approach revealed interesting drug protocols that remarkably reproduced the desired optimal effects simulated herein.
We illustrate the design of a generic new drug in Section~\ref{drugprotocols} by employing the equivalent dose and effects of docetaxel because this is the most common compound used for the chemotherapy of PCa, hence providing a reference for comparison during drug design. 
However, our framework can also accommodate specific drug doses and effects, e.g., by respectively choosing $d_c$ and calibrating $\beta_c$ with data in Eq.~\eqref{eq_u3}.
The design of therapeutic solutions could be further refined in multiple directions, e.g.,  by exploring alternative paradigms to model the dynamic effects of specific cytotoxic and antiangiogenic drugs  \cite{Gorelik2008,Hinow2009,Fister2003,Leszczynski2019,Yin2019,IrurzunArana2020}, by explicitly including drug toxicity (e.g., throughout the time integral of drug concentration \cite{Benzekry2013,Leszczynski2019}), by accounting for the synergistic action of drug combinations \cite{Meyer2019,IrurzunArana2020}, or by considering other forms of cancer treatment (e.g., the cytotoxic action of radiation therapy  \cite{Powathil2007,Corwin2013,Lima2017,Henares-Molina2017}).
Ultimately, the optimal pharmacodynamic properties calculated with our approach could assist in the development of new drug compounds that are more efficacious, show lower toxicities, and may even target specific cancer subtypes or patient-specific tumors \cite{Iyengar2012,Luepfert2005,Hussain2019,Shi2017}.
We also plan to explore the design of drug protocols for 3D PCa growth scenarios and considering longer simulation times, which are closer to the reality of experimental and clinical settings \cite{Kelly2012,Mottet2018,lorenzo2016tissue,lorenzo2017hierarchically,lorenzo2019computer,CGLMRR-1}.
Additionally, the acquisition of longitudinal PSA and imaging data during treatment would enable to recalibrate key parameters in our framework, for example: tumor proliferation and apoptosis in our PCa model, or drug effect rates $\beta_c$ and $\beta_a$ during drug protocol design. Although we have assumed these parameters to be constant, they are known to vary due to treatment action and the phenotypic evolution of the tumor \cite{Kim2005,Seruga2011,Hanahan2011,Gallaher2018,IrurzunArana2020}. 
Data-driven reparameterization would also enable us to update the optimal controls, adapt the design of personalized drug protocols accordingly, and early detect the emergence of chemoresistance patterns. 
In sum, we believe that our optimal control framework can provide a versatile and powerful approach to investigate patient-specific therapeutic strategies \emph{in silico} to succesfully treat advanced PCa.

\section*{Appendix A. Supplementary mathematical results.}\label{appA}
The notation we use in this appendix is somehow independent from the one in
the paper. We are going to check that if $u,v$ are two measurable functions
defined on $\Omega$ and $F:{\mathbb{R}} \to {\mathbb{R}}$ is a $C^1$
function, then there exists a measurable function $w:\Omega \to {\mathbb{R}}$
attaining intermediate values between the ones of $u$ and $v$ and such that 
\begin{equation}  \label{A1}
F (u(x)) - F (v(x)) = (u(x) - v(x)) F^{\prime }(w(x)) \quad \hbox{a.e. }
x\in \Omega.  \tag{A.1}
\end{equation}
The whole argument is due to Vittorino Pata, with many thanks from the
authors. \medskip

\noindent \textbf{Theorem A.1.} \textit{Let $f:{\mathbb{R}}\to{\mathbb{R}}$
be a continuous function. Then $f$ has a Borel measurable right inverse $k$
defined on (the interval) $\mathop{{\rm Im}} (f)$.}

\medskip

For the proof of Theorem~A.1 we need a lemma.

\medskip

\noindent \textbf{Lemma A.2.} \textit{Let $f:[a,b]\to{\mathbb{R}}$ be
continuous. Then the function $k:f([a,b])\to [a,b]$ given by} 
\begin{equation*}
k(t)=\inf\, \{y\in{\mathbb{R}} ; \ f(y)=t \}
\end{equation*}
\textit{is Borel measurable. Besides, for every $t\in f([a,b])$, we have that%
} 
\begin{equation*}
f(k(t))=t.
\end{equation*}

\medskip

\noindent \textbf{Proof. } The latter equality is obvious, since $f$ is
continuous, and so the infimum is actually a minimum. Let us prove the first
assertion. Let $t$ be fixed, such that $k(t)\neq a$. We will show that $k$
is either left or right continuous at $t$. This is enough to ensure that $k$
is Borel measurable. To this end, call $y=k(t)\in(a,b]$. By definition, $%
f(\tau)\neq y$ whenever $\tau<t$. Assume then that $f(\tau)<y$ for all $%
\tau<t$. We prove that $k$ is left continuous at $y$ (if instead $f(\tau)>y$
for all $\tau<t$ then $k$ is right continuous at $y$, the proof being the
same). By contradiction, if $k$ is not left continuous, there is $%
y_n\uparrow y$ such that $t_n=k(y_n)\not\to k(y)=t$. By compactness, there
is $\bar t<t$ such that $t_n\to\bar t$, up to a subsequence. But since $f$
is continuous, $y_n=f(t_n)\to f(\bar t)$, so that $f(\bar t)=y$. This
contradicts the fact that $t$ is the smallest element of $[a,b]$ for which
the equality $f(t)=y$ holds.\hfill $\square $

\medskip

\noindent \textbf{Proof of Theorem A.1. } For every $n\in{\mathbb{Z}}$, let 
\begin{equation*}
A_n^{\prime }=f([n,n+1)).
\end{equation*}
We now modify the sets $A_n^{\prime }$ in order to have them disjoint. Hence
we put 
\begin{equation*}
A_0=A^{\prime }_0,\quad A_1=A^{\prime }_1\setminus A_0,\quad
A_{-1}=A^{\prime }_{-1}\setminus(A_0\cup A_1), \quad A_{2}=A^{\prime
}_{2}\setminus(A_0\cup A_1\cup A_{-1}),
\end{equation*}
and so on. Possibly, some $A_n$ can be empty, and we simply ignore it.
Clearly, we have 
\begin{equation*}
\textstyle \bigcup_n A_n=\mathop{{\rm Im}}(f).
\end{equation*}
Let now $f_n$ be the restriction of $f$ on $[n,n+1]$, and let $%
k_n:f([n,n+1])\to [n,n+1]$ be given by 
\begin{equation*}
k_n(t)=\inf \{f_n^{-1}(t)\},
\end{equation*}
From the Lemma, $k_n$ is Borel measurable and $f_n(k_n(t))=t$ for all $t\in
f([n,n+1]).$ Finally, define the (Borel measurable) function 
\begin{equation*}
k(t)=k_n(t),\qquad t\in A_n.
\end{equation*}
Then, for $t\in A_n$, we have that $f(k(t))=f_n(k_n(t))=t$.\hfill $\square $

\medskip

\noindent \textbf{Remark A.3.} With inessential changes in the proof,
Theorem~A.1 still holds if $f$ is defined on a generic interval $I$.

\medskip

\noindent \textbf{Corollary A.4.} \textit{Given a domain $\Omega\subset{%
\mathbb{R}}^n$, let $\mu:\Omega\to {\mathbb{R}}$ be a measurable function
(finite everywhere), and let $f:{\mathbb{R}}\to{\mathbb{R}}$ be a continuous
function such that $\mathop{{\rm Im}}(f)\supset \mathop{{\rm Im}} (\mu)$.
Then there exists a measurable function $q:\Omega\to{\mathbb{R}}$ such that} 
\begin{equation*}
f(q(x))=\mu(x),\quad\hbox{for all } x\in\Omega.
\end{equation*}

\medskip

\noindent \textbf{Proof. } By Theorem~A.1, let $k$ be the Borel measurable
right inverse of $f$, and define 
\begin{equation*}
q(x)=k(\mu(x)),\quad x \in \Omega .
\end{equation*}
$q$ is Lebesgue measurable, being the composition of a Borel measurable
function and a Lebesgue measurable one.\hfill $\square $

\medskip

We finally conclude with the desired application involving the Lagrange mean
value theorem.

\medskip

\noindent \textbf{Corollary A.5.} \textit{Given a domain $\Omega \subset {%
\mathbb{R}}^n$, $n\in {\mathbb{N}}$, let $u,v : \Omega \to {\mathbb{R}}$ be
measurable functions (finite everywhere), and let $F : {\mathbb{R}} \to {%
\mathbb{R}}$ be a $C^1$ function. Then there exists a measurable function $%
w:\Omega \to {\mathbb{R}}$ such that} (\ref{A1}) \textit{holds.}\hfill $%
\square $

\medskip

\noindent \textbf{Proof. }We define 
\begin{equation*}
\mu(x)= 
\begin{cases}
\displaystyle \frac{F (u(x))- F (v(x))}{u(x) - v(x)} & \hbox{if } u(x)\not=
v(x), \\[0.3cm] 
F^{\prime }(u(x)) & \hbox{if } u(x)= v(x).%
\end{cases}%
\end{equation*}
Then the problem amounts to finding a measurable function $w$ such that 
\begin{equation*}
F^{\prime }(w(x)) = \mu(x), \quad x\in \Omega.
\end{equation*}
The Lagrange theorem implies that $\mathop{{\rm Im}} (F^{\prime }) \supseteq %
\mathop{{\rm Im}} (\mu).$ Hence, the conclusion follows from Corollary~A.4
by taking $w=q$.\hfill $\square $

\medskip

\newpage
\section*{Appendix B. Supplementary simulations.}\label{appB}

\renewcommand\thefigure{B.\arabic{figure}}    
\setcounter{figure}{0} 

\begin{figure}[!h]
\centerline{\includegraphics[width=\linewidth]{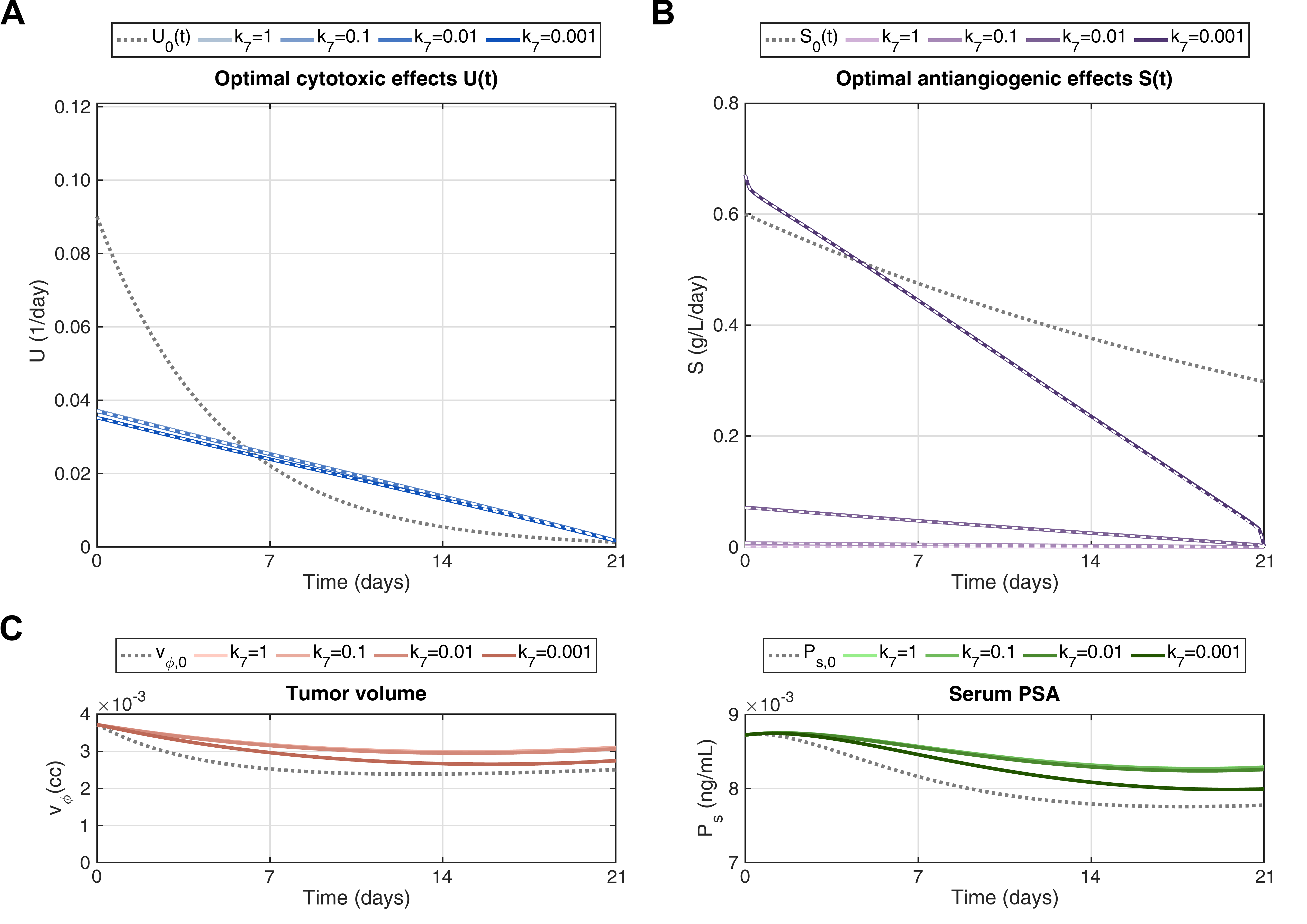}}
\vspace*{8pt}
\caption{Results of the optimal control problem using functional $J_1$ with $k_1=k_2=k_4=1$, $k_6=1$, and $k_7=1,10^{-1},10^{-2},10^{-3}$. (A) Optimal cytotoxic effects $U(t)$ obtained for each value of $k_7$ compared to standard docetaxel therapy ($U_0(t)$, gray dotted line) and the corresponding theoretical estimates calculated with Eq.~(\ref{40-0}) (overlapping white dashed lines). (B) Optimal antiangiogenic effects $S(t)$ obtained for each value of $k_7$ compared to standard bevacizumab therapy ($S_0(t)$, gray dotted line) and the corresponding theoretical estimates calculated with Eq.~(\ref{40-0}) (overlapping white dashed lines). (C) Evolution of tumor volume $v_\phi$ (left) and serum PSA $P_s$ (right) using the optimal $U(t)$ and $S(t)$ obtained via simulation for each value of $k_7$ compared to standard combined therapy, i.e., $v_{\phi,0}$ and $P_{s,0}$  respectively (gray dotted lines).  }
\end{figure}

\newpage

\begin{figure}[!h]
\centerline{\includegraphics[width=\linewidth]{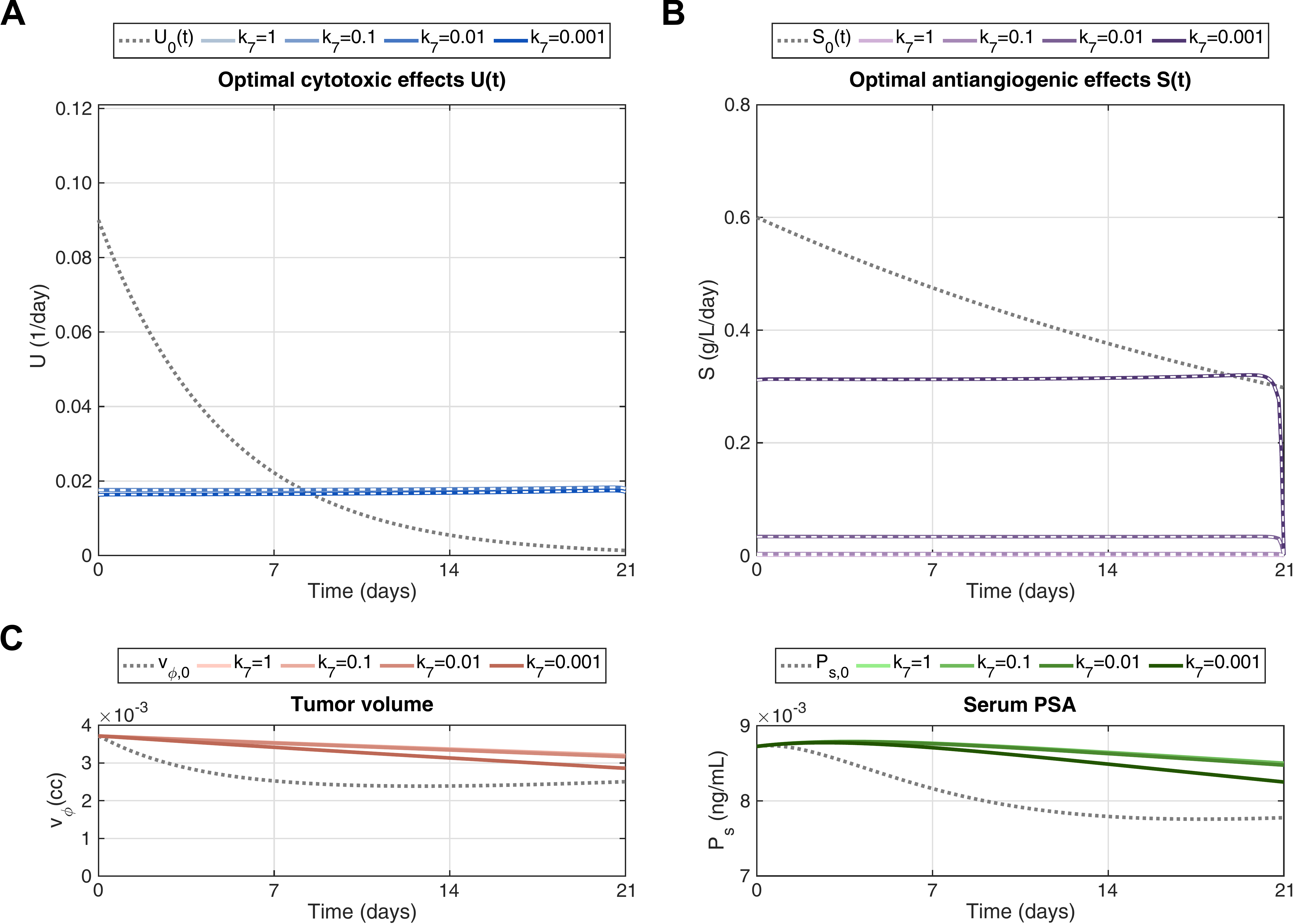}}
\vspace*{8pt}
\caption{Results of the optimal control problem using functional $J_2$ for $k_2=k_4=10$, $k_6=1$, and $k_7=1,10^{-1},10^{-2},10^{-3}$. (A) Optimal cytotoxic effects $U(t)$ obtained for each value of $k_7$ compared to standard docetaxel therapy ($U_0(t)$, gray dotted line) and the corresponding theoretical estimates calculated with Eq.~(\ref{40-0}) (overlapping white dashed lines). (B) Optimal antiangiogenic effects $S(t)$ obtained for each value of $k_7$ compared to standard bevacizumab therapy ($S_0(t)$, gray dotted line) and the corresponding theoretical estimates calculated with Eq.~(\ref{40-0}) (overlapping white dashed lines). (C) Evolution of tumor volume $v_\phi$ (left) and serum PSA $P_s$ (right) using the optimal $U(t)$ and $S(t)$ obtained via simulation for each value of $k_7$ compared to standard combined therapy, i.e., $v_{\phi,0}$ and $P_{s,0}$  respectively (gray dotted lines).   }
\end{figure}

\newpage

\begin{figure}[!h]
\centerline{\includegraphics[width=\linewidth]{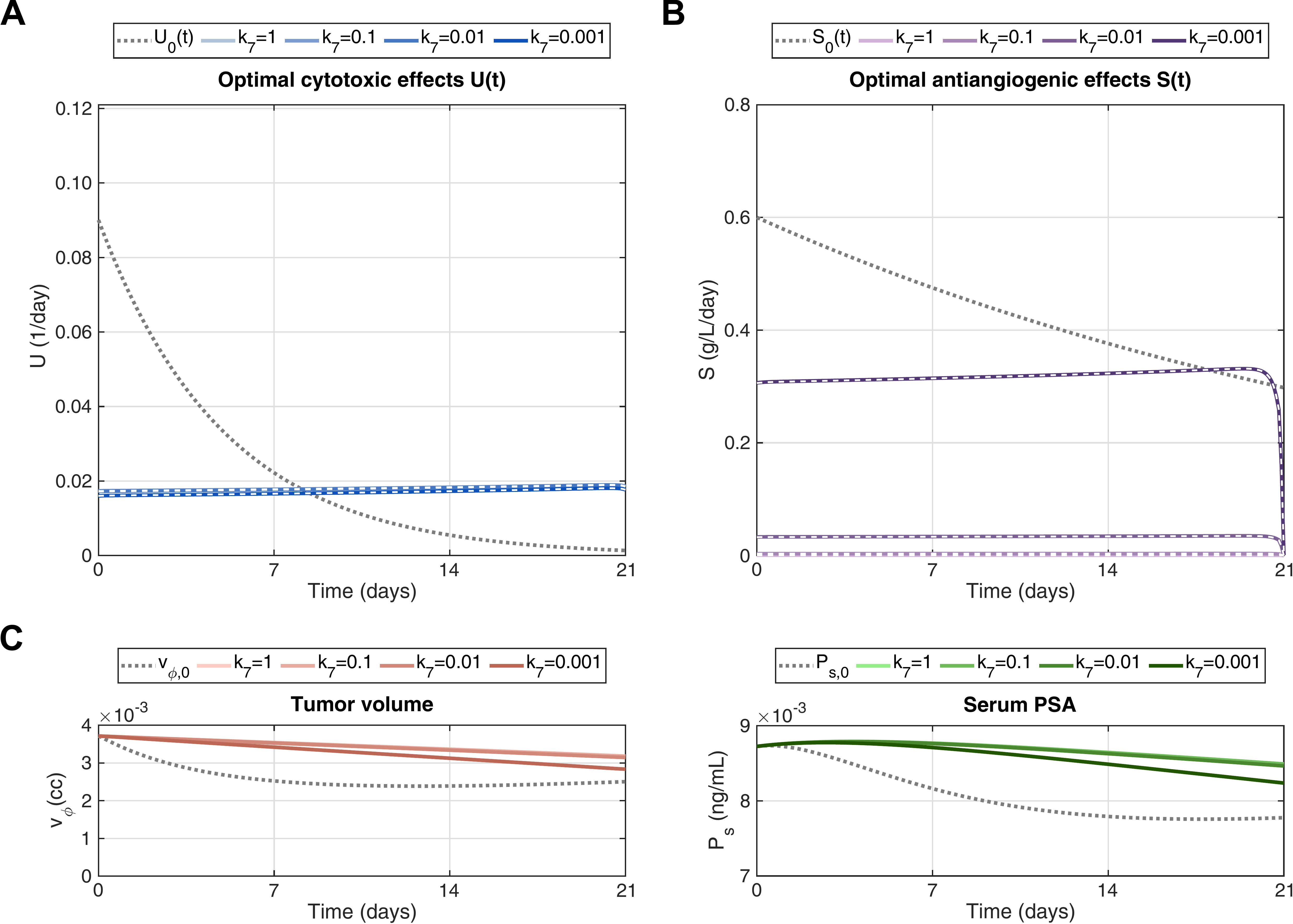}}
\vspace*{8pt}
\caption{Results of the optimal control problem using functional $J_3$ for $k_3=k_4=5$, $k_6=1$, and $k_7=1,10^{-1},10^{-2},10^{-3}$. (A) Optimal cytotoxic effects $U(t)$ obtained for each value of $k_7$ compared to standard docetaxel therapy ($U_0(t)$, gray dotted line) and the corresponding theoretical estimates calculated with Eq.~(\ref{40-0}) (overlapping white dashed lines). (B) Optimal antiangiogenic effects $S(t)$ obtained for each value of $k_7$ compared to standard bevacizumab therapy ($S_0(t)$, gray dotted line) and the corresponding theoretical estimates calculated with Eq.~(\ref{40-0}) (overlapping white dashed lines). (C) Evolution of tumor volume $v_\phi$ (left) and serum PSA $P_s$ (right) using the optimal $U(t)$ and $S(t)$ obtained via simulation for each value of $k_7$ compared to standard combined therapy, i.e., $v_{\phi,0}$ and $P_{s,0}$ respectively (gray dotted lines).   }
\end{figure}

\newpage

\section*{Acknowledgements}

The authors would like to express their gratitude to Viorel Barbu and Vittorino Pata for some fruitful discussions. 
This research was supported by the Italian Ministry of Education, University and Research (MIUR): Dipartimenti di Eccellenza Program (2018--2022) -- Dept. of Mathematics ``F. Casorati'', University of Pavia.
This research activity has been performed in the framework of the collaboration projects 
between the Italian CNR and the Romanian Academy: ``Control and stabilization problems for phase field and biological systems'' and ``Analysis and Optimization of mathematical models ranging from bio-medicine to engineering''.
The financial support of the project Fondazione Cariplo-Regione Lombardia MEGAs-TAR \textquotedblleft Matematica d'Eccellenza in biologia ed ingegneria come acceleratore di una nuova strateGia per l'ATtRattivit\`{a} dell'ateneo pavese\textquotedblright\ is gratefully acknowledged by ER. 
The paper also benefits from the support of the GNAMPA (Gruppo Nazionale per l'Analisi Matematica, la Probabilit\`{a} e le loro Applicazioni) of INdAM (Istituto Nazionale di Alta Matematica) for PC and ER. 
GL and AR have been partially supported by the MIUR-PRIN project XFAST-SIMS (no. 20173C478N). 
GL is also partially supported by a Peter O'Donnell Jr. Postdoctoral Fellowship from the Oden Institute for Computational Engineering and Sciences at The University of Texas at Austin.
The authors acknowledge the Rosen Center for Advanced Computing at Purdue University (USA) for providing HPC resources that contributed to the results presented in this paper. 

\bigskip

\begin{footnotesize}

\end{footnotesize}

\end{document}